\documentclass[a4paper,11pt]{article}
\usepackage{amsfonts}
\usepackage[latin1]{inputenc}
\usepackage{amssymb}
\usepackage{amsmath}
\usepackage[english]{babel}
\usepackage{fancyhdr}
\usepackage{cases}
\numberwithin{equation}{section}
\renewcommand{\baselinestretch}{1.40}
\newtheorem{thm}{Theorem}[section]
\newtheorem{lem}{Lemma}[section]
\newtheorem{prop}{Proposition}[section]

\def\no{\noindent}

\def \Max{\displaystyle\max}
\def\qed{ \hfill \vrule width.25cm height.25cm depth0cm\smallskip}
\input{tcilatex}
\begin{document}

\title{ Viscosity Solutions for a System of PDEs and\\ Optimal Switching  }

\author{Brahim EL ASRI \thanks{Universit\'e Ibn Zohr, Equipe. Aide à la decision,
ENSA, B.P.  1136, Agadir, Maroc. e-mail: b.elasri@uiz.ac.ma }\,\,\,
\, and \, Imade FAKHOURI
\thanks{Universit\'e Cadi Ayyad,  D\'ept.  de  Math\'ematiques, FSS,
B.P. 2390,  Marrakech, 40.000, Maroc. e-mail:
imadefakhouri@gmail.com.\newline
This author is supported by CNRST "Centre National pour la Recherche Scientifique et Technique", Rabat, Morocco.} }\maketitle

\begin{abstract}
In this paper, we study the existence and  uniqueness of  viscosity
solutions for a system of $m$ variational partial differential
inequalities with inter-connected obstacles. A particular case of
this system is the deterministic version of the Verification Theorem
of the Markovian optimal $m$-states optimal switching problem in
finite horizon. The switching cost functions are arbitrary and can
be positive or negative. This has an economic incentive in terms of
central valuation in cases where such organizations or state give
grants or financial assistance to power plants that promotes green
energy in their production activity or that uses less polluting
modes in their production. Our main tools is an approximation scheme
and the notion of systems of reflected backward stochastic
differential equations.
\end{abstract}

\noindent {$\mathbf{Keywords}$.} Real options, Backward stochastic differential equations, Snell envelope, Stopping times, Switching, Viscosity solution of PDEs, Variational inequalities

\medskip

\section{Introduction}
In this paper we consider the optimal m-states switching problem in finite horizon when the switching costs are
arbitrary and not necessarily  positive which is the novelty of this paper. This has an economical motivation in terms of firms valuation. \\
In order to introduce the problem let us deal with an example.
Assume a power plant which produces electricity and which has
several modes of production, therefore it is put in the
instantenuous most profitable one which is affected by the price
$(X_{t})_{t\geq 0}$ of electricity in the market that fluctuates in
reaction to many factors such as demand level, weather conditions,
unexpected outages etc, in addition to the fact that electricity is
non-storable so once produced, it should be immediately consumed.
Thus the manager of the plant aims at maximizing its global profit.
For this objective, she implements an optimal strategy which is a
pair of two sequences $(\tau_{n})_{n\geq 1}$ and
$(\epsilon_{n})_{n\geq 1}$describing respectively the optimal
successive switching times and modes. When the plant is in mode
$i\in \mathcal{I}$, it provides a profit $\psi_{i}(t,X_{t})dt$ which
depends on that mode. However this gain also incorporates a
switching cost $g_{ij}(t,X_{t})$, that could be positive or
negative, when switching the plant from the mode $i$ to another one.
This means that when $g_{ij}(t,X_{t})>0$, then switching is not free
and generates expenditures, on the other hand when $ g_ {ij} (t, X_
{t})\leq 0 $, it is the case when the state and environmental
organizations, provide grants and financial aid to power plants
that use the green energy in their production activities or methods of
cleaner production, which emit less carbon into the air.

The switching from one regime to another one is realized
sequentially at random times which are part of the decisions. So the
manager of the power plant faces two main issues:

$(i)$ when should she decide to switch the production from its
current mode to another one?

$(ii)$ to which mode the production has to be switched when the
decision of switching is made? \\
Optimal switching problems  were studied by several authors (see
e.g.
\cite{[BE],[BOU],[BO1],[BS],[CL],[DP],[DH],[DHP],[DZ2],[E1],[E2],[EH],[HJ],[LP],
[TY],[dz]} and the references therein). The motivations are mainly
related to decision making in the economic sphere. In order to
tackle those problems, authors use mainly two approaches. Either a
probabilistic one \cite{[DH], [DHP],[HJ],[HZ]} or an approach which uses
partial differential inequalities (PDIs for short)
\cite{[BO1],[CL],[DZ2],[EH],[HM],[TY],[dz]}.

In the finite horizon framework Djehiche \textit{et al.}
\cite{[DHP]} have studied the multi-modes switching problem when the
profit and the switching costs only depend on t, by using
probabilistic tools. They proved existence of a solution and found
an optimal strategy when the switching costs from state $i$ to state
$j$ is strictly non-negative ($g_{ij}(t)>\alpha >0$). The partial
differential equations approach (PDE in short) of this work has been carried out by
El Asri and Hamad\`{e}ne \cite{[EH]} when $g_{ij}(t,X_{t})>\alpha
>0$. They showed that when the price process $(X_{t}:t\geq 0)$ is
solution of a Markovian stochastic differential equation, then this
problem is associated to a system of variational inequalities with
interconnected obstacles for which they provided a solution in
viscosity sense. This solution turns out to be the value function of
the problem. Moreover the solution of the system is unique. In the
same spirit El Asri \cite{[E1]} studied the problem when
$g_{ij}(t,X_{t})\geq 0$, he showed the existence of the optimal
strategy and uniqueness of the solution in viscosity sense of the
problem. Nevertheless those papers \cite{[E1]} and \cite{[EH]}
suffer from two facts: (i) the switching cost functions $g_{ij}$ are
non negative ; (ii) in the markovian case the optimal strategy
satisfy  $\forall n\geq 1,$ $p[\tau_n<T]< \frac{C}{n}$ in
\cite{[EH]}, this property is the basis of the proof of existence of
a solution in viscosity sense, but in the general case where the
costs of switching could be negative this property is not satisfied.
In the PDEs approach we also mention the recent result of Hamad\`ene
and Morlais \cite{[HM]} that deals with existence and uniqueness,
when the switching cost functions are positive and arbitrary, in
viscosity sense of a solution for a $m$ system of variational PDIs
with interconnected obstacles which is the deterministic version of
the Verification Theorem of the Markovian optimal switching problem.

The novelty of this paper lies in the fact that we investigate the
solution to the optimal multiple switching problem when the
switching costs could be positive or negative, using probabilistic
tools as the Snell envelope of processes, backward stochastic
differential equations (BSDEs for short), and partial differential
equations approach.

 We prove existence and uniqueness of the vector
of value functions and  provide a characterization of an optimal
strategy of this problem when the payoff rates $\psi _{i}$ and the
switching costs $g_{ij}$ (positive or negative) are adapted only to
the filtration generated by a Brownian motion. Later on, in the
markovian framework, we show that the value function of the problem
is associated to an uplet of deterministic functions
$(v_{1},\dots,v_{m})$ which is the unique solution of the following
system of PDIs:
\begin{equation}
\left\{
\begin{array}{l}
\min\left\{v_{i}(t,x)-\max\limits_{j\in \mathcal{I}^{-i}}\left\{
-g_{ij}(t,x)+v_{j}(t,x)\right\} ,-\partial_{t}v_{i}(t,x)-\mathcal{A}v_{i}(t,x)-\psi _{i}(t,x)\right\} =0\\
\forall \,\,(t,x)\in \lbrack 0,T]\times \mathbb{R}^{k},\,\,i\in
\mathcal{I}=\{1,...,m\},\text{\thinspace\ and \thinspace\
}v_{i}(T,x)=0,
\end{array}
\right.
\end{equation}
where $\mathcal{A}$ an operator associated with a diffusion process
and $\mathcal{I}^{-i}:=\mathcal{I}\setminus \{i\}$. It turns out
that this system is the deterministic version of the Verification
Theorem of the optimal multi-modes switching problem in finite
horizon.

This paper is organized as follows: \\
In Section 2, we formulate the problem and give the related
definitions. In Section 3, we shall introduce the optimal switching
problem under consideration and give its probabilistic Verification
Theorem. It is expressed by means of the Snell envelope of
processes. Then we introduce the approximating scheme which enables
us to construct a solution for the Verification Theorem. %Moreover we
%give some properties of that solution, especially the dynamic
%programming principle .
Section 4 is devoted to the connection
between the optimal switching problem, the Verification Theorem and
the associated system of PDIs. This connection is made through BSDEs
with one reflecting obstacle in the Markovian case. Further we show
existence and continuity of a solution for the system of PDIs.
Finally, in Section 5, we show that the solution of PDIs is unique
in the class of continuous functions which satisfy a polynomial
growth condition.
\section{ Formulation of the problem and assumptions}
\subsection{Setting of the problem}
The finite horizon multiple switching problem can be formulated as
follows. Let $\mathcal{I}$ be the set of all possible activity modes
of the production of a power plant. A management strategy of the
plant consists, on the one hand, of the choice of a sequence of
nondecreasing stopping times $(\tau _{n})_{n\geq 1}$ (i.e. $\tau
_{n}\leq \tau _{n+1}$, $\tau _{0}=0$ and $\tau _{n}\rightarrow T $
when $ n\rightarrow +\infty$) where the manager decides to switch
the activity from its current mode to another one. On the other
hand, it consists of the choice of the mode $\xi_{n}$, which is an
$\mathcal{F}_{\tau_{n}}$-measurable random variable taking values in
$\mathcal{I}$, to which the production is switched at $\tau _{n}$
from its current mode. Therefore the admissible management
strategies of the plant are the pairs $(\delta ,\xi ):=((\tau
_{n})_{n\geq 1},(\xi_{n})_{n\geq 1})$ and the set of these
strategies is denoted by $\mathcal{D}$.

Let $X:=(X_{t})_{0\leq t\leq T}$ be an adapted continuous %$\mathcal{P}$-measurable,
$\mathbb{R}^{k}$-valued  stochastic process, which stands
for the market price of $k$ factors which determine the market price
of the commodity. Assuming that the production activity is in mode 1
at the initial time $t=0$, let $(u_{t})_{0\leq t\leq T}$ denote the
indicator of the production activity's mode at time $t\in \lbrack
0,T]$:
\begin{equation}
u_{t}=\mathbf{1}_{[0,\tau _{1}]}(t)+\sum_{n\geq
1}\xi_{n}\mathbf{1}_{(\tau _{n},\tau _{n+1}]}(t).
\end{equation}
Then for any $t\leq T$, the state of the whole economic system
related to the project at time $t$ is represented by the vector
\begin{equation}
(t,X_{t},u_{t})\in \lbrack 0,T]\times \mathbb{R}^{k}\times
\mathcal{I}.
\end{equation}
Finally, let $\psi _{i}(t,X_{t})$ be the instantaneous profit when
the system is in state $(t,X_{t},i)$, and for $i,j\in \mathcal{I}$
$\ i\neq j$, let $g_{ij}(t,X_{t})$ denote the switching cost of the
production at time $t$ from current mode $i$ to another mode $j$.
Then if the plant is run under the strategy $(\delta ,\xi
)=((\tau_{n})_{n\geq 1},(\xi _{n})_{n\geq 1})$ the expected total
profit is given by:
\begin{equation*}
J(\delta ,\xi
)=E\left[\int_{0}^{T}\psi_{u_{s}}(s,X_{s})ds-\sum_{n\geq
1}g_{u_{\tau_{n-1}}u_{\tau_{n}}}(\tau _{n},X_{\tau
_{n}})\mathbf{1}_{[\tau_{n}<T]}\right].
\end{equation*}
Therefore the problem we are interested in is to find an optimal
strategy $i.e.$ a strategy $(\delta^{\ast },\xi^{\ast })$ such that
$J(\delta^{\ast },\xi ^{\ast })\geq J(\delta ,\xi )$ for any
$(\delta ,\xi )\in \mathcal{D}$.

We now consider the following system of $m$ variational inequalities
with inter-connected obstacles: $\forall \,\,i\in \mathcal{I}$
\begin{equation}\label{pdi}
\left\{
\begin{array}{l}
\min \left\{ v_{i}(t,x)-\max\limits_{j\in \mathcal{I}^{-i}}\left\{-g_{ij}(t,x)+v_{j}(t,x)\right\} ,-\partial_{t}v_{i}(t,x)-\mathcal{A}v_{i}(t,x)-\psi _{i}(t,x)\right\} =0, \\
v_{i}(T,x)=0,
\end{array}
\right.
\end{equation}
where $\mathcal{A}$ is given by:
\begin{equation}
\label{derivgen} \mathcal{A}=\frac{1}{2}\sum_{i,j=1}^{m}(\sigma
\sigma^{\ast })_{ij}(t,x)\frac{\partial^{2}}{\partial x_{i} \partial
x_{j}}+\sum_{i=1}^{m}b_{i}(t,x) \frac{\partial }{\partial x_{i}}\,;
\end{equation}
hereafter the superscript $(^{\ast })$ stands for the transpose,
$Tr$ is the trace operator and finally $\left\langle x,y\right\rangle $ is the inner product of $x,y\in \mathbb{R}^{k}$.

The main objective of this paper is to focus on the existence and
uniqueness of the solution in viscosity sense of (\ref{pdi}). This
system is the deterministic version of the optimal $m$-states
switching problem when we assume that the market price process $X$
of the commodity is an It\^{o} diffusion.

Recall the notion of viscosity solution of the system (\ref{pdi}).
\begin{definition}
Let $(v_1,\ldots ,v_m)$ be a $m$-uplet of continuous functions
defined on $[0,T]\times \mathbb{R}^{k}$, $\mathbb{R}$-valued and
such that $v_{i}(T,x)=0$ for any $x\in \mathbb{R}^{k}$ and $i\in
\mathcal{I}$. The $m$-uplet $(v_{1},\ldots ,v_{m})$ is called:
\begin{itemize}
\item[$(i)$] a viscosity supersolution (resp. subsolution) of the system (\ref{pdi})
if for each fixed $i\in \mathcal{I}$, for any $(t_{0},x_{0})\in
[0,T]\times \mathbb{R}^{k}$ and any function $\varphi _{i}\in
\mathcal{C}^{1,2}([0,T]\times \mathbb{R}^{k})$ such that $\varphi
_{i}(t_{0},x_{0})=v_{i}(t_{0},x_{0})$ and $(t_{0},x_{0})$ is a local
maximum of $\varphi _{i}-v_{i}$ (resp. minimum), we have:
\begin{eqnarray}
&&\min \bigg\{v_{i}(t_{0},x_{0})-\max_{j\in\mathcal{I}^{-i}}\left(-g_{ij}(t_{0},x_{0})+v_{j}(t_{0},x_{0})\right), \nonumber\\ &&\qquad \qquad \qquad-\partial_{t}\varphi_{i}(t_{0},x_{0})-\mathcal{A}\varphi_{i}(t_{0},x_{0})-\psi_{i}(t_{0},x_{0})
\bigg\} \geq 0 \,\,(\mbox{resp.}\leq 0);
\end{eqnarray}

\item[$(ii)$] a viscosity solution if it is both a viscosity supersolution and subsolution. $\qed$
\end{itemize}
\end{definition}

There is an equivalent formulation of this definition (see
e.g.\cite{[CIL]}) which we give because it will be useful later. So
firstly we define the notions of superjet and subjet of a continuous
function $v$.

\begin{definition}
Let $v\in \mathcal{C}((0,T)\times \mathbb{R}^{k})$, $(t,x)$ an
element of $(0,T)\times \mathbb{R}^{k}$ and finally $\mathbf{S}_{k}$
the set of $k\times k$ symmetric matrices. We denote by
$J^{2,+}v(t,x)$ (resp.$J^{2,-}v(t,x)$), the superjets (resp. the
subjets) of $v$ at $(t,x)$, the set of triples $(p,q,X)\in
\mathbb{R}\times \mathbb{R}^{k}\times \mathbf{S}_{k}$ such that:
\begin{eqnarray*}
v(s,y) &\leq &v(t,x)+p(s-t)+\langle q,y-x\rangle +\frac{1}{2}\langle X(y-x),y-x\rangle +o(|s-t|+|y-x|^{2}) \\
&& \\
\bigg(resp.\ v(s,y) &\geq &v(t,x)+p(s-t)+\langle q,y-x\rangle
+\frac{1}{2}\langle X(y-x),y-x\rangle +o(|s-t|+|y-x|^{2})\bigg).
\end{eqnarray*}
\end{definition}

Note that if $\varphi-v$ has a local maximum (resp. minimum) at
$(t,x)$, then we obviously have:
\begin{equation*}
\bigg(D_{t}\varphi
(t,x),D_{x}\varphi(t,x),D_{xx}^{2}\varphi(t,x)\bigg) \in
J^{2,-}v(t,x)\,\,\,(\mbox{resp.}\ J^{2,+}v(t,x)).\qed
\end{equation*}

We now give an equivalent definition of a viscosity solution of the
parabolic system with inter-connected obstacles.

\begin{definition}
Let $(v_{1},\ldots ,v_{m})$ be a $m$-uplet of continuous functions
defined on $[0,T]\times \mathbb{R}^k$, $\mathbb{R}$-valued and such
that
$(v_{1},\ldots ,v_{m})(T,x)=0$ for any $x\in \mathbb{R}^{k}$. The $m$-uplet $(v_{1},\ldots ,v_{m})$ is called a viscosity supersolution (resp. subsolution) of (\ref{pdi})
if for any $i\in \mathcal{I}$, $(t,x)\in (0,T)\times \mathbb{R}^k$
and $(p,q,X)\in J^{2,-} v_i (t,x)$ (resp. $J^{2,+} v_i (t,x)$),
\begin{equation*}
min \left\{v_i(t,x)-\max\limits_{j\in\mathcal{I}^{-i}}(-g_{ij}(t,x)+ v_{j}(t,x)),-p -\frac{1}{2}Tr[\sigma^{*} X \sigma] -\langle b,q\rangle-\psi_{i}(t,x)\right\}\geq 0 \,\,(resp.\leq 0).
\end{equation*}
It is called a viscosity solution if it is both a viscosity subsolution and supersolution. $\qed$
\end{definition}
\subsection{Assumptions}

Throughout this paper $T$ (resp. $k,\, d$) is a fixed real
(resp.integers) positive numbers. Let $(\Omega, \mathcal{F}, \mathbb{P})$ be
a fixed probability space on which is defined a standard
$d$-dimensional Brownian motion $B=(B_{t})_{0\leq t\leq T}$ whose
natural filtration is $(\mathcal{F}_{t}^{0}:=\sigma \{B_{s}, s\leq
t\})_{0\leq t\leq T}$. Let $\mathbf{F}=(\mathcal{F}_{t})_{0\leq
t\leq T}$ be the completed filtration of
$(\mathcal{F}_{t}^{0})_{0\leq t\leq T}$ with the $\mathbb{P}$-null sets of
$\mathcal{F}$. \\

Furthermore, let:

- $\mathcal{P}$ be the $\sigma$-algebra on $[0,T]\times \Omega$ of
$\mathbf{F}$-progressively measurable sets;

- $\mathcal{M}^{2,k}$ be the set of $\mathcal{P}$-measurable and
$\mathbb{R}^{k}$-valued processes $w=(w_{t})_{t\leq T}$ such that
$E[\int_{0}^{T}|w_{s}|^{2}ds]<\infty$ and $\mathcal{S}^2$ be the set
of $\mathcal{P}$-measurable, continuous processes ${w}=({w}_t)_{t\leq
T}$ such that $E[\sup_{t\leq T}|{w}_{t}|^2]<\infty$;

- for any stopping time $\tau \in [0,T]$, $\mathcal{T}_\tau$ denotes
the set of all stopping times $\theta$ such that $\tau \leq \theta
\leq T$. \medskip

- $\Pi$ be the class of functions with polynomial growth, defined as follows:
\begin{eqnarray*}
{\Pi}&:=&\bigg\{\varphi:(t,x)\in [0,T]\times \mathbb{R}^{k}\rightarrow
\varphi(t,x) \in \mathbb{R}, \quad such \ that \\&& \quad|\varphi(t,x)|\leq C(1+|x|^\gamma) \mbox{ for
some non negative real constants } C \mbox{ and }\gamma\bigg\}.
\end{eqnarray*}

We now make the following assumptions on the data:\\
\noindent {\bf [H1]}: $b:[0,T]\times \mathbb{R}^{k}\rightarrow \mathbb{R}^{k}$ and $\sigma :[0,T]\times \mathbb{R}^{k}\rightarrow \mathbb{R}^{k\times d}$ are two continuous functions for which there exists a constant $C>0$ such that for any $t\in \lbrack 0,T]$ and $x,x^{\prime }\in \mathbb{R}^{k}$
\begin{eqnarray}\label{b}
&&|\sigma (t,x)|+|b(t,x)|\leq C (1+|x|) \nonumber  \ \ \mbox{and} \\&& \nonumber \\
&&|\sigma (t,x)-\sigma (t,x^{\prime })|+|b(t,x)-b(t,x^{\prime})|\leq C|x-x^{\prime }|.
\end{eqnarray}
Throughout this paper we assume that assumption [H1] holds. \medskip

\no {\bf [H2]}: For $i\in \mathcal{I}$,
$\psi_{i}:[0,T]\times\mathbb{R}^{k}\rightarrow \mathbb{R}$ is
continuous and belongs to $\Pi$.\medskip

\no {\bf [H3]}: For any $i,j\in \mathcal{I}$ and $(t,x)\in [0,T]\times \mathbb{R}^{k}$:
\begin{itemize}
\item[$(i)$] $g_{ij}:[0,T]\times
\mathbb{R}^{k} \rightarrow \mathbb{R}$ is jointly continuous in $(t,x)$, could be positive or negative and belongs to $\Pi$. Moreover as a convention we assume that $g_{ii}(t,x)=0$.
\item[$(ii)$] For any sequence of indices $i_{1},\ldots,i_{k}\in \mathcal{I}$ such that $i_{1}=i_{k}$ and $card\{i_1,...,i_k\}=k-1$ we have:
\begin{equation}
\label{nofreeloop}
g_{i_{1}i_{2}}(t,x)+g_{i_{2}i_{3}}(t,x)+\ldots+g_{i_{k-1}i_{k}}(t,x)+g_{i_{k}i_{1}}(t,x)>0.
\end{equation}
\item [$(iii)$] For any $i,j \in \mathcal{I}$, if $g_{ij}$ is non positive, we assume that
\begin{equation}\label{cost-ter}
g_{ij}(T,x)=0, \quad j\neq i.
\end{equation}
\end{itemize}
\section{The Verification Theorem and existence of the  processes $Y^{i},i=1,\ldots,m$}

Note that in order that the quantity $J(\delta ,\xi )$ makes sense, we assume throughout this paper that, for any $i,j\in \mathcal{I}$ the processes
$(\psi _{i}(t,X_{t}))_{t\leq T}$ and $(g_{ij}(t,X_{t}))_{t\leq T}$ belong to $\mathcal{M}^{2,1}$ and $\mathcal{S}^{2}$ respectively. There is
one to one correspondence between the pairs $(\delta ,\xi )$ and the pairs $(\delta ,u)$. Therefore throughout this paper one refers indifferently to
$(\delta ,\xi )$ or $(\delta ,u)$.

%\begin{defn}\label{admissible}
%The admissible management strategies are the pairs $(\delta,\xi)=((\tau_n)_{n\geq 1},(\xi_n)_{n\geq 1})$ that verify the assumption (\ref{cumulative-switching-costs}).\\
%The set of admissible management strategies is denoted $\mathcal{D}$ and defined:
% \begin{equation*}
% \mathcal{D}=\left\{(\delta ,\xi )=
% ((\tau_n)_{n\geq 1},(\xi_n)_{n\geq 1})\ s.t \quad E\Big|\sum_{n\geq1}g_{\xi_{n-1}\xi_{n}}(\tau _{n},X_{\tau_{n}})\mathbf{1}_{[\tau_{n}\leq T]}\Big|^2<\infty \right\}.
% \end{equation*}
%
%\end{defn}

\subsection{The Verification Theorem}

To tackle the problem described above Djehiche \textit{et al.}
\cite{[DHP]} have introduced a Verification Theorem which is expressed
by means of Snell envelope of processes. The Snell envelope of a
stochastic process $(\eta _{t})_{0\leq t\leq T}$ of $\mathcal{S}^{2}$
(with a possible positive jump at $T$)
is the lowest supermartingale $R(\eta ):=(R(\eta )_{t})_{0\leq t\leq T}$ of $
\mathcal{S}^{2}$ such that for any $t\leq T$, $R(\eta )_{t}\geq \eta
_{t}$. It has the following expression:
\begin{equation*}
\forall \;t\leq T,\ R(\eta )_{t}=\limfunc{ess\ sup}_{\tau \in \mathcal{T}
_{t}}E[\eta _{\tau }|\mathcal{F}_{t}]\mbox{ and satisfies }R(\eta
)_{T}=\eta _{T}.
\end{equation*}
For more details we refer to \cite{[CK], [Elka]}. \medskip

The Verification Theorem for the $m$-states optimal switching
problem is the following:
\begin{thm}
\label{thmverif} (Verification Theorem) \\
Assume that for any $i,j \in \mathcal{I}$ the following hold:
\begin{itemize}
\item [(a)] $\psi_i$ satisfies (H2);

\item [(b)] $g_{ij}$ satisfies (H3)-(i) and (ii);

\item [(c)] there exist $m$ processes $
(Y^{i}:=(Y_{t}^{i})_{0\leq t\leq T},i=1,\ldots ,m)$ of
$\mathcal{S}^{2}$ such that:
\begin{eqnarray}\label{yi}
\forall \;t\leq T,\ \ && Y_{t}^{i}=\limfunc{ess\ sup}_{\tau \geq t}E\left[\int_{t}^{\tau}\psi_{i}(s,X_{s})ds+\max\limits_{j\in \mathcal{I}
^{-i}}(-g_{ij}(\tau ,X_{\tau })+Y_{\tau }^{j})\mathbf{1}_{[\tau <T]}|\mathcal{F}_{t}\right], \nonumber \\&&
Y_{T}^{i}=0.
\end{eqnarray}
\end{itemize}
%\begin{equation}\label{yi}
%\left\{
%\begin{array}{l}
%\forall \;t\leq T,\ \ Y_{t}^{i}=\limfunc{ess\ sup}_{\tau \geq t}E\left[\int_{t}^{\tau}\psi_{i}(s,X_{s})ds+\max\limits_{j\in \mathcal{I}
%^{-i}}(-g_{ij}(\tau ,X_{\tau })+Y_{\tau }^{j})\mathbf{1}_{[\tau <T]}|\mathcal{F}_{t}\right],\\
%Y_{T}^{i}=0.
%\end{array}
%\right.
%\end{equation}
Then:

\begin{itemize}
\item[(i)] $Y_{0}^{1}=\sup\limits_{(\delta ,\xi )\in \mathcal{D}}J(\delta,u).$

\item[(ii)] Define the sequence of $\mathbf{F}$-stopping times $\delta^{\ast }=(\tau _{n}^{\ast })_{n\geq 1}$ as follows :
\begin{equation*}
\begin{array}{lll}
\tau _{1}^{\ast } & = & \inf \left\{ s\geq 0,\quad Y_{s}^{1}=\max\limits_{j\in {\mathcal{\ I}^{-1}}}(-g_{1j}(s,X_{s})+Y_{s}^{j})\right\} \wedge T.\\
\mbox{For} \quad n\geq 2, \\
\tau_{n}^{\ast} & = & \inf \left\{s\geq \tau_{n-1}^{\ast},\quad Y_{s}^{u_{\tau_{n-1}^{\ast}}}=\max\limits_{k\in \mathcal{I}\backslash \{u_{\tau_{n-1}^{\ast }}\}}(-g_{u_{\tau_{n-1}^{\ast}}k}(s,X_{s})+Y_{s}^{k})\right\} \wedge T,\ \\
\end{array}
\end{equation*}
where:

\begin{itemize}
\item[$\bullet $] $u_{\tau _{1}^{\ast }}=\sum\limits_{j\in \mathcal{\ I}}j%
\mathbf{1}_{\{\max\limits_{k\in \mathcal{\ I}^{-1}}(-g_{1k}(\tau
_{1}^{\ast },X_{\tau _{1}^{\ast }})+Y_{\tau _{1}^{\ast
}}^{k})=-g_{1j}(\tau _{1}^{\ast },X_{\tau _{1}^{\ast }})+Y_{\tau
_{1}^{\ast }}^{j}\}};$

\item[$\bullet $] for any $n\geq 1$ and $t\geq \tau _{n}^{\ast },$ $
Y_{t}^{u_{\tau _{n}^{\ast }}}=\sum\limits_{j\in \mathcal{I}}\mathbf{1}
_{[u_{\tau _{n}^{\ast }}=j]}Y_{t}^{j}$;

\item[$\bullet $] for any $n\geq 2,\,\,u_{\tau _{n}^{\ast }}=l$ on the set
\begin{equation*}
\left\{ \max\limits_{k\in \mathcal{I}\backslash \{{u_{\tau _{n-1}^{\ast }}}\}}(-g_{u_{\tau _{n-1}^{\ast }}k}(\tau _{n}^{\ast },X_{\tau_{n}^{\ast }})+Y_{\tau _{n}^{\ast }}^{k})=-g_{u_{\tau _{n-1}^{\ast}l}}(\tau _{n}^{\ast },X_{\tau _{n}^{\ast }})+Y_{\tau _{n}^{\ast}}^{l}\right\}
\end{equation*}
with $g_{u_{\tau _{n-1}^{\ast }k}}(\tau _{n}^{\ast },X_{\tau_{n}^{\ast }})=\sum\limits_{j\in \mathcal{I}}\mathbf{1}_{[u_{\tau_{n-1}^{\ast}}
=j]}g_{jk}(\tau _{n}^{\ast },X_{\tau _{n}^{\ast }})$ \ and \ $\mathcal{I}\backslash \{u_{\tau _{n-1}^{\ast }}\}=\sum\limits_{j\in \mathcal{I}}\mathbf{1}_{[u_{\tau _{n-1}^{\ast }}=j]}\mathcal{I}^{-j}$.\newline
 Then the strategy $(\delta ^{\ast },u^{\ast})$ is optimal and admissible i.e. it satisfies:
$$-\infty<E\left[ \sum\limits_{k\geq 1}g_{u_{\tau_{k-1}^{\ast }}u_{\tau_{k}^{\ast }}}({\tau _{k}^{\ast }},X_{{\tau _{k}^{\ast}}})\mathbf{1}_{[\tau _{k}^{\ast }<T]}\right] <+\infty \qquad \mbox{and} \qquad \mathbb{P}[\tau_{n}^{\ast}<T \ , \ \forall \ n\geq 0]=0.$$
\end{itemize}
\end{itemize}
\end{thm}
$Proof .$  The proofs of $\textbf{(i)}$ and $\textbf{(ii)}$ are omitted since they are similar to the ones in \cite{[DHP], [E1]}, except for the property $\mathbb{P}[\tau_{n}^{\ast}<T \ , \ \forall \ n\geq 0]=0$ that we are going to give its proof below.

 Let us show now that the optimal strategy $(\delta ^{\ast },u^{\ast })=(\tau_{n}^{\ast},u_{\tau _{n}^{\ast }})_{n\geq 0}$ is admissible, that is we should prove that $\mathbb{P}[\tau_{n}^{\ast}<T \ , \ \forall \ n\geq 0]=0$. Let's assume the contrary, i.e. $\mathbb{P}[\tau_{n}^{\ast}<T \ , \ \forall \ n\geq 0]>0$ and show by contradiction that is impossible. \\
Let $\tau_{0}^{\ast}=0$ and $u_{0}^{\ast}=i$. Thanks to the definition of $\tau_{n}^{\ast}$ , we have: $$\mathbb{P}[Y_{\tau_{n+1}^{\ast}}^{u_{\tau _{n}^{\ast }}}=-g_{u_{\tau_{n}^{\ast}}u_{\tau_{n+1}^{\ast}}}(\tau_{n+1}^{\ast},X_{\tau_{n+1}^{\ast}})+Y_{\tau_{n+1}^{\ast}}^{u_{\tau_{n+1}^{\ast}}}, \ u_{\tau_{n+1}^{\ast}} \in \mathcal{I}^{-u_{\tau_{n}^{\ast}}} ,\ \forall n \geq 1]>0.$$
Since $\mathcal{I}$ is finite then there is a state $i_{0} \in \mathcal{I}$ and a loop $i_{0},i_{1},\ldots,i_{k},i_{0}$ of elements of $\mathcal{I}$ such that $card\{i_{0},i_{1},\ldots,i_{k}\}=k+1$ and:
$$\mathbb{P}[Y_{\tau_{n+1}^{\ast}}^{i_{l}}=-g_{i_{l}i_{l+1}}(\tau_{n+1}^{\ast},X_{\tau_{n+1}^{\ast}})+Y_{\tau_{n+1}^{\ast}}^{i_{l+1}}, \ l=0,\ldots,k \ ,(i_{k+1}=i_{0}), \ \forall \ n\geq 1]>0.$$
Therefore taking the limit w.r.t. $n$ to obtain:
$$\mathbb{P}[Y_{\tau}^{i_{l}}=-g_{i_{l}i_{l+1}}(\tau,X_{\tau})+Y_{\tau}^{i_{l+1}},
 \ l=0,\ldots,k \ ,(i_{k+1}=i_{0})]>0,$$
where $\tau:=\lim\limits_{n\rightarrow \infty} \tau_{n}^{\ast}$. But
this implies that:
$$\mathbb{P}[g_{i_{0}i_{1}}(\tau,X_{\tau})+ \ldots +g_{i_{k}i_{0}}(\tau,X_{\tau})=0]>0,$$
which contradicts (\ref{nofreeloop}). Therefore  $\mathbb{P}[\tau_{n}^{\ast}<T \ , \ \forall \ n\geq 0]=0$. \\ Then the optimal strategy is admissible.\ $\qed$
\begin{remark}
The condition $\mathbb{P}[\tau_{n}^{\ast}<T \ , \ \forall \ n\geq 0]=0$ means that the sequence $(\tau_{n}^{\ast}(\omega))_{n\geq 0}$ is stationary, moreover in economic terms it signifies that the manager is allowed to make only a finite number of decisions during time interval $[0,T]$, otherwise the switching costs would be infinite, and then $J(\delta ,\xi )$ would go to $0$.
\end{remark}
%%%%%%%%%%%%%%%

\subsection{Existence of the  processes $Y^{i},i=1,\ldots,m$}
we will now establish existence of the processes $Y^{1},\ldots ,Y^{m}$. They will be obtained as a limit of a
sequence of processes $(Y^{1,n},\ldots ,Y^{m,n})$ defined recursively by means of the Snell envelope notion as follows: \\
For $i\in \mathcal{I}$, we set, for any $0\leq t \leq T,$
\begin{equation}
Y_{t}^{i,0}=E\left[
\int_{t}^{T}\psi_{i}(s,X_{s})ds|\mathcal{F}_{t}\right] ,\,\,0\leq
t\leq T, \label{y0}
\end{equation}
and for $n\geq 1$,
\begin{equation}
\label{yin} Y_{t}^{i,n}=\limfunc{ess\ sup}_{\tau \in \mathcal{T}_t}E\left[
\int_{t}^{\tau}\psi
_{i}(s,X_{s})ds+\max\limits_{k\in\mathcal{I}^{-i}}(-g_{ik}(\tau,X_{\tau})+Y_{\tau}^{k,n-1})\mathbf{1}_{[\tau
<T]}|\mathcal{F}_{t}\right] ,\,\,0\leq t\leq T.
\end{equation}
Next we will give some useful properties of
$Y^{1,n},\ldots,Y^{m,n}.$
\begin{lem} Assume that for any $i,j \in \mathcal{I}$:
\begin{itemize}
\item [(i)] $\psi_i$ satisfies (H2);

\item [(ii)] $g_{ij}$ satisfies (H3)-(i) and (iii).
\end{itemize}
 Then, for any $n\geq 0$ the processes $Y^{1,n},\ldots,Y^{m,n}$ are continuous and belong to $\mathcal{S}^{2}$.
\end{lem}
$Proof.$ Let us show by induction that for any $n\geq 0$ and every $i \in \mathcal{I}$,\ the $Y^{i,n}$'s are continuous and belong to $\mathcal{S}^{2}$. \\
For $n=0$ the property holds true since we can write $Y^{i,0}$ as the sum of a continuous process and a martingale w.r.t to the Brownian filtration which is continuous, therefore $Y^{i,0}$ is continuous and since the process $(\psi_{i}(s,X_{s}))_{0 \leq s \leq T} $ belongs to $\mathcal{M}^{2}$ then by using Doob's inequality we obtain that  $Y^{i,0}$ belong to $\mathcal{S}^{2}$. \\
Suppose now that the property is satisfied for some $n$: \\
For every $i \in \mathcal{I}$ and up to a term, $Y^{i,n+1}$ is the
Snell envelope of the process: \\
$\left(\int_{0}^{t}\psi_{i}(s,X_{s})ds
+\max\limits_{k\in\mathcal{I}^{-i}}(-g_{ik}(t,X_{t})+Y_{t}^{k,n})\mathbf{1}_{[t<T]}\right)_{0 \leq t \leq T}$ and verifies $Y_{T}^{i,n+1}=0$. \\
The process
$\max\limits_{k\in\mathcal{I}^{-i}}(-g_{ik}(t,X_{t})+Y_{t}^{k,n})\mathbf{1}_{[t<T]}$
is continuous on $[0,T)$ thanks to the continuity of $Y_{t}^{k,n}$,
and at $T$ we have two cases: for every $i,k \in \mathcal{I}$
\begin{itemize}
\item[(a)] if $g_{ik}(T,X_{T})$ is positive then $ \max\limits_{k\in\mathcal{I}^{-i}}(-g_{ik}(t,X_{t})+Y_{t}^{k,n})|_{t=T}<0$. Thus the process  $\max\limits_{k\in\mathcal{I}^{-i}}(-g_{ik}(t,X_{t})+Y_{t}^{k,n})$ is continuous on $[0,T)$ and has a positive jump at T since $Y_{T}^{i,n+1}=0$. Then we deduce by Proposition 2 (iii) in \cite{[DHP]} that $Y^{i,n+1}$ is continuous on $[0,T]$.
\item[(b)]if $g_{ik}(T,X_{T})$ is negative then $\max\limits_{k\in\mathcal{I}^{-i}}(-g_{ik}(t,X_{t})+Y_{t}^{k,n})|_{t=T}=0$ since we have, by assumption (H3)-(iii), that $g_{ik}(T,X_{T})=0$ for $k\neq i$. Thus $Y^{i,n+1}$ is continuous on $[0,T]$.
\end{itemize}
Therefore we deduce that $Y^{i,n+1}$ is continuous on $[0,T]$ and belongs to $\mathcal{S}^{2}$. \\
This shows that for any $n\geq 0$ and every $i \in \mathcal{I}$, the
$Y^{i,n}$'s are continuous and belong to $\mathcal{S}^{2}$.$\qed$\\
%%%%%%%%%%%%%%%%%%%%%%%%%%%%%%%%%%%%%%%%%%%%%%%%%%%
In the Proposition \ref{prop3.1} below, we will show that the sequence of processes $(Y^{1,n},\ldots ,Y^{m,n})$ converge increasingly and pointwisely P-a.s. in $\mathcal{M}^{2,1}$. But to do so we will need an additional assumption on the negative switching costs:\\
\no {\bf [H4]}: $\forall n\geq 1$, for any $\tau_n\in \mathcal{T}_0$, $\xi_n\in\mathcal{I}$ and $x\in\mathbb{R}^{k}$ there exists some real $K>0$ such that:
\begin{equation}
Card\Big\{ 0\leq\tau_n\leq T ,\quad n\geq 1 \quad \mbox{such that} \quad g_{\xi_{n-1}\xi_n}(\tau_n,x)<0\Big\}\leq K.
\end{equation}
This assumption means that the number of stopping times where we can get negative switching costs is bounded which means in particular that the number of negative switching costs is also bounded. Actually economically speaking this assumption is realistic since the negative switching costs could be seen as a kind of reward in the form of grants or financial aid given, for example, to power plants using green energy as explained in the introduction. Notice that this financial aid could not be infinite (it should be bounded), otherwise that would mean that we could earn money by switching modes of production over and over but this will lead to an infinite value.\qed \\
We will discuss in the remark below from a mathematical point of view the reason why we impose this assumption.
\begin{remark}\label{H4-math-meaning}
We need to control the following quantity $$E\left[\sum_{j=1}^{n}(-g_{u_{\tau_{j-1}u_{\tau_{j}}}}
(\tau_{j},X_{\tau_{j}}))\mathbf{1}_{[\tau_{j}<T]}\right].$$
In fact, we only need to control that quantity over the negative switching costs, since we have the following
\begin{eqnarray*}
&&E\left[\sum_{j=1}^{n}\Big(-g_{u_{\tau_{j-1}u_{\tau_{j}}}}
(\tau_{j},X_{\tau_{j}})\Big)\mathbf{1}_{[\tau_{j}<T]}\right]\\
&\leq&E\left[\sum_{j=1}^{n}\Big(-g_{u_{\tau_{j-1}u_{\tau_{j}}}}
(\tau_{j},X_{\tau_{j}})\Big)\mathbf{1}_{[g_{u_{\tau_{j-1}u_{\tau_{j}}}}(\tau_{j},X_{\tau_{j}})\geq0]}\mathbf{1}_{[\tau_{j}<T]}\right]
\\&&+E\left[\sum_{j=1}^{n}\Big(-g_{u_{\tau_{j-1}u_{\tau_{j}}}}
(\tau_{j},X_{\tau_{j}})\Big)\mathbf{1}_{[g_{u_{\tau_{j-1}u_{\tau_{j}}}}(\tau_{j},X_{\tau_{j}})<0]}\mathbf{1}_{[\tau_{j}<T]}\right]
\\&\le&E\left[\sum_{j=1}^{n}\Big(-g_{u_{\tau_{j-1}u_{\tau_{j}}}}
(\tau_{j},X_{\tau_{j}})\Big)\mathbf{1}_{[g_{u_{\tau_{j-1}u_{\tau_{j}}}}(\tau_{j},X_{\tau_{j}})<0]}\mathbf{1}_{[\tau_{j}<T]}\right]
\end{eqnarray*}
Now applying assumption (H4) on the last term of the previous inequality, we obtain that
 \begin{equation}\label{H4-neg-sw-cost}
E\left[\sum_{j=1}^{n}(-g_{u_{\tau_{j-1}u_{\tau_{j}}}}
(\tau_{j},X_{\tau_{j}}))\mathbf{1}_{[\tau_{j}<T]}\right]\leq K E\left[\max_{l,j \in \mathcal{I}; l\neq j}\Big\{\sup_{s\le T}\Big(-g_{lj}(s,X_s)\Big)\Big\}\right].
\end{equation}

Notice that without assumption (H4) we can only get that
\begin{equation*}
E\left[\sum_{j=1}^{n}(-g_{u_{\tau_{j-1}u_{\tau_{j}}}}
(\tau_{j},X_{\tau_{j}}))\mathbf{1}_{[\tau_{j}<T]}\right]<\infty.
\end{equation*}
\end{remark}
%%%%%%%%%%%%%%%%%%%%%%%%%%%%%%%%%%%
\begin{prop}
\label{prop3.1}
 Assume that (H2), (H3) and (H4) are fullfiled. Then for any $i \in \mathcal{I}$, the sequence $(Y^{i,n})_{n\geq0}$ converges increasingly and pointwisely P-a.s. for any $0\leq t \leq T$ and in $\mathcal{M}^{2,1}$ to c\`{a}dl\`{a}g processes $\tilde{Y}^{i}$. \\ Moreover these limit processes $\tilde{Y}^{i}=(\tilde{Y}^{i}_{t})_{0 \leq t \leq T}$, $i=1,\ldots,m$, satisfy the following:
    \begin{itemize}
    \item[$(a)$] $ E\left[\sup\limits_{0\leq t \leq T}|\tilde{Y}^{i}_{t}|^{2}\right] < +\infty $ , $i \in \mathcal{I} $.
    \item[$(b)$] For any $0\leq t\leq T$ we have:
    \begin{equation}
    \label{ytilde}
    \tilde{Y}^{i}_{t}=\limfunc{ess\ sup}_{\tau \in \mathcal{T}_t} E\left[\int_{t}^{\tau}\psi_{i}(s,X_{s})ds+
    \max\limits_{k\in\mathcal{I}^{-i}}(-g_{ik}(\tau,X_{\tau})+\tilde{Y}^{k}_{\tau})\mathbf{1}_{[\tau <
    T]}|\mathcal{F}_{t}\right].
    \end{equation}
    \end{itemize}
\end{prop}
$Proof. $

 Let us now set $ \mathcal{D}_{t}^{i,n}=\{u \in \mathcal{D} \quad such\ that \quad u_{0}=i,\ \tau_{1} \geq t \quad and \quad \tau_{n+1}=T \} $. \\
Using the same arguments as the ones of the Verification Theorem,
the following characterization of the processes $Y^{i,n}$ holds
true:
\begin{equation*}
Y_{t}^{i,n}=\limfunc{ess\ sup}_{u\in\mathcal{D}_{t}^{i,n}}E\left[\int_{t}^{T}\psi_{u_{s}}(s,X_{s})ds-\sum_{j=1}^{n}g_{u_{\tau_{j-1}u_{\tau_{j}}}}(\tau_{j},X_{\tau_{j}})
\mathbf{1}_{[\tau_{j}<T]}|\mathcal{F}_{t}\right]
\end{equation*}
Since $\mathcal{D}_{t}^{i,n}\subset\mathcal{D}_{t}^{i,n+1}$, we have
P-a.s. for all $t\in[0,T]$, $Y_{t}^{i,n}\leq Y_{t}^{i,n+1}$ thanks
to the continuity of $Y^{i,n}$. Moreover we have:
\begin{eqnarray}
Y_{t}^{i,n} &=&
\limfunc{ess\ sup}_{u\in\mathcal{D}_{t}^{i,n}}E\left[\int_{t}^{T}\psi_{u_{s}}(s,X_{s})ds-\sum_{j=1}^{n}g_{u_{\tau_{j-1}u_{\tau_{j}}}}
(\tau_{j},X_{\tau_{j}})\mathbf{1}_{[\tau_{j}<T]}|\mathcal{F}_{t}\right]  \nonumber \\
&\leq&\limfunc{ess\ sup}_{u\in\mathcal{D}_{t}^{i,n}}E\left[\int_{t}^{T}\psi_{u_{s}}(s,X_{s})ds|\mathcal{F}_{t}\right]+\limfunc{ess\ sup}_{u\in\mathcal{D}_{t}^{i,n}}E\left[\sum_{j=1}^{n}\Big(-g_{u_{\tau_{j-1}u_{\tau_{j}}}}
(\tau_{j},X_{\tau_{j}})\Big)\mathbf{1}_{[\tau_{j}<T]}|\mathcal{F}_{t}\right]\nonumber\\
&\leq& E\left[\int_{t}^{T}\limfunc{max}_{i \in \mathcal{I}}
|\psi_{i}(s,X_{s})|ds|\mathcal{F}_{t}\right] + \limfunc{ess\ sup}_{u\in\mathcal{D}_{t}^{i,n}}E\left[\sum_{j=1}^{n}\Big(-g_{u_{\tau_{j-1}u_{\tau_{j}}}}
(\tau_{j},X_{\tau_{j}})\Big)\mathbf{1}_{[\tau_{j}<T]}|\mathcal{F}_{t}\right].
\end{eqnarray}
Now using assumption (H4) and Remark \ref{H4-math-meaning} (inequality (\ref{H4-neg-sw-cost})), we get that
\begin{eqnarray*}
Y_{t}^{i,n} &\leq& E\left(\int_{t}^{T}\limfunc{max}_{i \in
\mathcal{I}}|\psi_{i}(s,X_{s})|ds|\mathcal{F}_{t}\right) + K
E\left[\limfunc{max}_{l,j \in \mathcal{I};\ j\neq
l}\Big\{\limfunc{sup}_{s\leq
T}\Big(-g_{lj}(s,X_{s})\Big)\Big\}|\mathcal{F}_{t}\right].
\end{eqnarray*}
Therefore, for every $i \in \mathcal{I}$ the sequence
$(Y^{i,n})_{n\geq 0}$ satisfies
\begin{equation}\label{increasing-inequa}
Y_{t}^{i,n} \leq Y_{t}^{i,n+1} \leq E\left(\int_{t}^{T}\limfunc{max}_{i \in
\mathcal{I}}|\psi_{i}(s,X_{s})|ds|\mathcal{F}_{t}\right) + K
E\left[\limfunc{max}_{l,j \in \mathcal{I};\ j\neq
l}\Big\{\limfunc{sup}_{s\leq
T}\Big(-g_{lj}(s,X_{s})\Big)\Big\}|\mathcal{F}_{t}\right]; \ \forall t\in[0,T].
\end{equation}
Now since $\psi_i$ and $g_{ij}$ belong to $\mathcal{M}^{2,1}$ and $\mathcal{S}^{2}$ respectively.
Then $(Y^{i,n})_{n\geq 0}$ converges to some limit
$\tilde{Y}^{i}_{t}:=\lim\limits_{n \rightarrow +\infty}Y^{i,n}_{t}$
that satisfies: $\forall i \in \mathcal{I}$
\begin{equation}\label{Y-tilde}
Y^{i,0}_{t} \leq \tilde{Y}^{i}_{t} \leq
E\left(\int_{t}^{T}\limfunc{max}_{i \in
\mathcal{I}}|\psi_{i}(s,X_{s})|ds|\mathcal{F}_{t}\right) + K
E\left[\limfunc{max}_{l,j \in \mathcal{I};\ j\neq
l}\Big\{\limfunc{sup}_{s\leq
T}\Big(-g_{lj}(s,X_{s})\Big)\Big\}|\mathcal{F}_{t}\right]; \quad \forall t\in[0,T].
\end{equation}
Next, using (\ref{Y-tilde}), the fact that $\psi_i$ and $g_{ij}$ belong to $\mathcal{M}^{2,1}$ and $\mathcal{S}^{2}$ respectively, and
Doob's Maximal Inequality yield, for each $i\in \mathcal{I}$
$$E\Big(\limfunc{sup}_{t\leq T}|\tilde{Y}^{i}_{t}|^{2}\Big)< +\infty.$$
By the Lebesgue Dominated Convergence Theorem, the sequence $(Y^{i,n})_{n\geq 0}$ also converges to $\tilde{Y}^{i}$ in $\mathcal{M}^{2,1}$. \\
Let us now show that $\tilde{Y}^{i}$ is c\`{a}dl\`{a}g. Actually, for each
 $i\in \mathcal{I}$ and $n\geq 1$, by (\ref{yin}) the process $\left(Y^{i,n}_{t}+\int_{0}^{t}\psi(s,X_{s})ds\right)_{0\leq t\leq T}$ is
 a continuous supermartingale. Hence its limit process \\
 $\left(\tilde{Y}^{i}_{t}+\int_{0}^{t}\psi(s,X_{s})ds\right)_{0\leq
t\leq T}$ is c\`{a}dl\`{a}g as a limit of increasing sequence of
continuous
 supermartingales (see Dellacherie and Meyer \cite{[DM]} p.86). Therefore $\tilde{Y}^{i}$ is c\`{a}dl\`{a}g. \\
Finally, since the c\`{a}dl\`{a}g processes
$\tilde{Y}^{1},\ldots,\tilde{Y}^{m}$ are limits of the sequence of
increasing continuous processes $Y^{i,n},i\in \mathcal{I}$ that
satisfy (\ref{yin}), then by Snell envelope properties, the
processes $\tilde{Y}^{1},\ldots,\tilde{Y}^{m}$ satisfy the following: for any $0\leq t\leq T$,\quad $i=1,\ldots,m$
\begin{equation*}
    \tilde{Y}^{i}_{t}=\limfunc{ess\ sup}_{\tau \in \mathcal{T}_t} E\left[\int_{t}^{\tau}\psi_{i}(s,X_{s})ds+
    \max\limits_{k\in\mathcal{I}^{-i}}(-g_{ik}(\tau,X_{\tau})+\tilde{Y}^{k}_{\tau})\mathbf{1}_{[\tau <
    T]}|\mathcal{F}_{t}\right],
    \end{equation*}
    which is the desired result.\qed \medskip

We will now prove that the processes
$\tilde{Y}^{1},\ldots,\tilde{Y}^{m}$ are continuous and satisfy the
Verification Theorem (Theorem \ref{thmverif}).
\begin{thm} Assume that (H2), (H3) and (H4) hold. Then
the limit processes $\tilde{Y}^{1},\ldots,\tilde{Y}^{m}$ satisfy the
Verification Theorem (Theorem \ref{thmverif}).
\end{thm}
$Proof.$  Recall from Proposition \ref{prop3.1} that the processes $\tilde{Y}^{1},\ldots,\tilde{Y}^{m}$ are c\`{a}dl\`{a}g and uniformly ${L}^{2}$--integrable and satisfy (\ref{ytilde}). It remains to prove that they are continuous. \\
Indeed, note that, for $i\in \mathcal{I}$, the process $(\tilde{Y}^{i}_{t}+\int^{t}_{0}\psi_{i}(s)ds)_{0\leq t\leq T} $ is the Snell envelope of \\
$$\left(\int^{t}_{0}\psi_{i}(s,X_{s})ds+\max\limits_{k\in\mathcal{I}^{-i}}(-g_{ik}(t,X_{t})+\tilde{Y}^{k}_{t})\mathbf{1}_{[t<T]}\right)_{0\leq t\leq T}$$
since the processes $(\int^{t}_{0}\psi_{i}(s,X_{s})ds)_{0\leq t\leq T}$ are continuous. Therefore from the property of the jumps of the Snell envelope (see \cite{[DHP]}; Proposition 2 (ii)),
when there is a (necessarily negative) jump of $\tilde{Y}^{i}_{t}$ at $t$, there is a jump, at the same time t, of the process $\left(\max\limits_{k\in\mathcal{I}^{-i}}(-g_{ik}(t,X_{t})+\tilde{Y}^{k}_{t})\right)_{0\leq t < T}$. Since $g_{ij}$ are continuous, there is $j\in \mathcal{I}^{-i}$ such that $\Delta_{t}\tilde{Y}^{j}<0$, and $\tilde{Y}^{i}_{t-}=-g_{ij}(t,X_{t})+\tilde{Y}^{j}_{t-}$. \\
Suppose now there is an index $i_{1} \in \mathcal{I}$ for which
there exists $t\in [0,T]$ such that $\Delta_{t}\tilde{Y}^{i_{1}}<0$.
This implies that there exists another index $i_{2} \in
\mathcal{I}^{-i_{1}}$ such that $\Delta_{t}\tilde{Y}^{i_{2}}<0$ and
$\tilde{Y}^{i_{1}}_{t-}=-g_{i_{1}i_{2}}(t,X_{t})+\tilde{Y}^{i_{2}}_{t-}$. But given $i_{2}$, there exists an index $i_{3}\in \mathcal{I}^{-i_{2}}$ such that $\Delta_{t}\tilde{Y}^{i_{3}}<0$ and $\tilde{Y}^{i_{2}}_{t-}=-g_{i_{2}i_{3}}(t,X_{t})+\tilde{Y}^{i_{3}}_{t-}$. Repeating this argument many times, we get a sequence of indices $i_{1},...,i_{j},... \in \mathcal{I}$ that have the property that  $i_{k}\in \mathcal{I}^{-i_{k-1}}$, $\Delta_{t}\tilde{Y}^{i_{k}}<0$ and $\tilde{Y}^{i_{k-1}}_{t-}=-g_{i_{k-1}i_{k}}(t,X_{t})+\tilde{Y}^{i_{k}}_{t-}$. \\
Since $\mathcal{I}$ is finite, there exist two indices $q<r$ such
that $i_{q}=i_{r}$ and $i_{q},i_{q+1},...,i_{r-1}$ are mutually
different. It follows that
\begin{eqnarray*}
\tilde{Y}^{i_{q}}_{t-}&=&-g_{i_{q}i_{q+1}}(t,X_{t})+\tilde{Y}^{i_{q+1}}_{t-} \\
&=&-g_{i_{q}i_{q+1}}(t,X_{t})-g_{i_{q+1}i_{q+2}}(t,X_{t})+\tilde{Y}^{i_{q+2}}_{t-} \\
&=&... \\
&=&-g_{i_{q}i_{q+1}}(t,X_{t})-...-g_{i_{r-1}i_{r}}(t,X_{t})+\tilde{Y}^{i_{r}}_{t-}.
\end{eqnarray*}
As $i_{q}=i_{r}$ we get
$$-g_{i_{q}i_{q+1}}(t,X_{t})-...-g_{i_{r-1}i_{r}}(t,X_{t})=0,$$ which
contradicts assumption (\ref{nofreeloop}).
Therefore there is no $i \in \mathcal{I}$ for which there is a $t\in [0,T]$ such that $\Delta_{t}\tilde{Y}^{i}<0$. This means that the processes  $\tilde{Y}^{1},\ldots,\tilde{Y}^{m}$ are continuous. Since they satisfy (\ref{ytilde}), then by uniqueness $Y^{i}=\tilde{Y}^{i}$ for any $i\in \mathcal{I}$. \\
Thus the Verification Theorem (Theorem \ref{thmverif}) is satisfied by $Y^1,\ldots,Y^m$. \qed
%$\tilde{Y}^{1},\ldots,\tilde{Y}^{m}$.$\qed $

\section{Existence of a solution for the system of variational inequalities}
In this section we will address the question of existence of a solution for the system of variational inequalities (\ref{pdi}). But first let's make the link between those solutions and BSDEs with one reflecting barrier in the Markovian framework.
\subsection{Connection with BSDEs with  one reflecting barrier}
Let $(t,x)\in [0,T]\times \mathbb{R}^{k}$ and let $(X^{tx}_s)_{s\leq T}$ be the solution of the following standard SDE:
\begin{equation}
\label{sde}
dX^{tx}_{s}=b(s,X_{s}^{tx})\,ds+\sigma(s,X_{s}^{tx})\,dB_{s} \mbox{ for }\ t\leq
s\leq T\mbox{ and }\ X_{s}^{tx}=x \mbox{ for }\ s\leq t\end{equation}where the functions $b$ and $\sigma$ are the ones of (\ref{b}). These properties of $\sigma$ and $b$ imply in particular that the process $(X^{tx}_s)_{0\le s\leq T}$ solution of the standard SDE (\ref{sde})  exists and is unique, for any $t\in [0, T]$ and $x\in \mathbb{R}^{k}$.

The operator $\mathcal{A}$  that appears in (\ref{derivgen}) is the infinitesimal generator associated with $X^{tx}$. In the
following result we collect some properties of $X^{tx}$.

\begin{prop}
(see e.g. \cite{[RY]}) The process $X^{tx}$ satisfies the following estimates:
\begin{itemize}
\item [$(i)$] For any $q\geq 2$, there exists a constant $C$ such that
\begin{equation}\label{estimat1}
E[\sup_{0\le s\leq T}|X^{tx}_s|^q]\leq C(1+|x|^q).
\end{equation}
\item[$(ii)$] There exists a constant $C$ such that for any $t,t'\in [0,T]$ and $x,x'\in \mathbb{R}^{k}$,
\begin{equation}\label{estimat2}
E[\sup_{0\le s\leq T}|X_{s}^{tx}-X_{s}^{t'x'}|^2]\leq C(1+|x|^2)(|x-x'|^2+|t-t'|). \qed
\end{equation}
\end{itemize}
\end{prop}

We are going now to introduce the notion of a BSDE with one reflecting barrier introduced in \cite{[EKal]}. This notion will
allow us to make the connection between the variational inequalities system  and the $m$-states optimal switching problem
described in the previous section.
\medskip

So let us introduce the deterministic functions $f:[0,T]\times \mathbb{R}^{k+1+d}\rightarrow \mathbb{R}$, \\ $h:[0,T]\times \mathbb{R}^{k}\rightarrow \mathbb{R}$ and $g:\mathbb{R}^{k}\rightarrow \mathbb{R}$ continuous, of polynomial growth and such that $h(x,T)\leq g(x)$. Moreover we assume that for any $(t,x)\in [0,T]\times \mathbb{R}^{k}$, the mapping $(y,z)\in \mathbb{R}^{1+d}\mapsto f(t,x,y,z)$ is uniformly Lipschitz. Then we have the following
result related to BSDEs with one reflecting barrier:
\begin{thm}\label{yitx}(\cite{[EKal]}, Th. 5.2 and 8.5) For any $(t,x)\in [0,T]\times \mathbb{R}^{k}$, there exits a unique triple of processes $(Y^{tx},Z^{tx},K^{tx})$ such that:
\begin{equation}
\left\{
\begin{array}{l}
Y^{tx}, K^{tx}\in {\cal S}^2 \mbox{ and }Z^{tx}\in \mathcal{M}^{2,d};\,K^{tx} \mbox{ is  non-decreasing and }K^{tx}_0=0,\\
Y^{tx}_{s}=g(X^{tx}_{T})+\int_{s}^{T}f(r,X_r^{tx},Y^{tx}_{r},Z^{tx}_{r})dr-\int_{s}^{T}Z^{tx}_{r}dB_{r}+
K_{T}^{tx}-K^{tx}_{s}, \,\, s\leq T,\\
Y^{tx}_{s}\geq h(s,X_{s}^{tx}),\, \forall s\leq T\mbox{ and }
\int_{0}^{T}(Y^{tx}_r-h(r,X^{tx}_{r}))dK^{tx}_r=0.
\end{array}
\right.
\end{equation}
 Moreover the following characterization of $Y^{tx}$ as a Snell envelope holds true:
\begin{equation}
\label{snellenv}
\forall s\leq T,\,\,Y^{tx}_{s}=\limfunc{ess\ sup}_{\tau\in \mathcal{T}_{s}}E\left[\int_{t}^{\tau} f(r,X_{r}^{tx},Y^{tx}_{r},Z^{tx}_{r})dr+h(\tau,X^{tx}_{\tau})\mathbf{1}_{[\tau <T]}+g(X^{tx}_{T})\mathbf{1}_{[\tau=T]}|\mathcal{F}_{s}\right].
\end{equation}

On the other hand there exists a deterministic continuous with
polynomial growth function $u: [0,T]\times \mathbb{R}^{k}
\rightarrow \mathbb{R}$ such that:$$\forall s\in [t,T],
\quad Y^{tx}_s=u(s,X^{tx}_s).$$ Moreover the function $u$ is the unique
viscosity solution in the class of continuous functions with
polynomial growth of the following PDE with obstacle:
$$
\left\{
\begin{array}{l}
\min\{u(t,x)- h(t,x), -\partial_{t}u(t,x)-
\mathcal{A}u(t,x)-f(t,x,u(t,x),\sigma (t,x)^*\nabla u(t,x))\}=0,
\\u(T,x)=g(x).\qed
\end{array}
\right.
$$
\end{thm}
\subsection{Existence of a solution for the system of variational inequalities}

Let $(Y^{1,tx}_s,...,Y^{m,tx}_s)_{0\leq s\leq T}$ be the processes which satisfy the Verification Theorem \ref{thmverif}  in the case when the process $X\equiv X^{tx}$. Therefore using the characterization (\ref{snellenv}),
there exist processes $K^{i,tx}$ and $Z^{i,tx}$, $i\in \mathcal{I}$, such that the triples ($Y^{i,tx},Z^{i,tx},K^{i,tx})$ are unique solutions of the following reflected BSDEs: for any $i=1,...,m$ we have,
\begin{equation}
\left\{
\begin{array}{l}
Y^{i,tx}, K^{i,tx}\in {\cal S}^2 \mbox{ and }Z^{i,tx}\in \mathcal{M}^{2,d};\,K^{i,tx} \mbox{ is  non-decreasing and }K^{i,tx}_0=0,\\
Y^{i,tx}_{s}=\int_{s}^{T}\psi_i(r,X_{r}^{tx})du-\int_{s}^{T}Z^{i,tx}_{r}dB_{r}+K_{T}^{i,tx}-K^{i,tx}_{s}, \,\,\, 0\leq s\leq T,\,\,Y^{i,tx}_{T}=0,\\
Y^{i,tx}_{s}\geq \max\limits_{j\in \mathcal{I}^{-i}}(-g_{ij}(s,X_{s}^{tx})+Y^{j,tx}_{s}),\,\,0\leq s\leq T,\\
\int_{0}^{T}(Y^{i,tx}_{r}-\max\limits_{j\in \mathcal{I}^{-i}}(-g_{ij}(r,X_{r}^{tx})+Y^{j,tx}_{r}))dK^{i,tx}_{r}=0.
\end{array}
\right.
\end{equation}
Moreover we have the following result.
\begin{prop}\label{viscosity}
Assume that (H2), (H3) and (H4) hold. Then there are deterministic functions $v^{1},...,v^{m}$ $:[0,T]\times \mathbb{R}^{k} \rightarrow \mathbb{R}$ such that:
$$\forall (t,x)\in [0,T]\times \mathbb{R}^{k},\ \forall s\in [t,T],\ Y_s^{i,tx}=v_i(s,X_{s}^{tx}), \,\,i=1,...,m.$$
Moreover the functions $v_i$, $i=1,...,m,$ are of polynomial growth.
\end{prop}
$Proof.$ For $n\geq 0$ let
$(Y^{n,1,tx}_{s},...,Y^{n,m,tx}_{s})_{0\leq s\leq T}$ be the
processes constructed in (\ref{y0})-(\ref{yin}). Therefore using an
induction argument and Theorem \ref{yitx} there exist deterministic
continuous with polynomial growth functions $v_{i,n}$ $(i=1,...,m)$
such that for any $(t,x)\in [0,T]\times \mathbb{R}^{k}$, $\forall
s\in [t,T]$, $Y^{i,n,tx}_{s}=v_{i,n}(s,X_{s}^{tx})$.
Next since we have that
\begin{equation*}
Y_{t}^{i,n,tx} =
\limfunc{ess\ sup}_{u\in\mathcal{D}_{t}^{i,n}}E\left[\int_{t}^{T}\psi_{u_{s}}(s,X_{s}^{tx})ds+\sum_{j=1}^{n}(-g_{u_{\tau_{j-1}u_{\tau_{j}}}}
(\tau_{j},X_{\tau_{j}}^{tx}))\mathbf{1}_{[\tau_{j}<T]}|\mathcal{F}_{t}\right],
\end{equation*}
then following the same reasoning as in Remark \ref{H4-math-meaning} and Proposition \ref{prop3.1} (inequality \ref{increasing-inequa}) we obtain that
\begin{equation*}
Y^{i,n,tx}_t\le Y^{i,n+1,tx}_t\leq E\left[\int_{t}^{T}\limfunc{max}_{i \in \mathcal{I}}|\psi_{i}(s,X_{s}^{tx})|ds|\mathcal{F}_{t}\right]
 + KE\left[\limfunc{max}_{l,j \in \mathcal{I};\ j\neq l}\left(\limfunc{sup}_{s\leq T}(-g_{lj}(s,X_{s}^{tx}))\right)|\mathcal{F}_{t}\right].
\end{equation*}
Next, since $Y^{i,n,tx}_{t}$ and $Y^{i,n+1,tx}_{t}$ are deterministics, then taking expectations in the last inequality leads to
\begin{equation*}
Y^{i,n,tx}_t\le Y^{i,n+1,tx}_t\leq E\left[\int_{t}^{T}\limfunc{max}_{i \in \mathcal{I}}|\psi_{i}(s,X_{s}^{tx})|ds\right]
 + KE\left[\limfunc{max}_{l,j \in \mathcal{I};\ j\neq l}\left(\limfunc{sup}_{s\leq T}(-g_{lj}(s,X_{s}^{tx}))\right)\right].
\end{equation*}
Therefore combining the polynomial growth of $\psi_i$ and $g_{ij}$ and estimate (\ref{estimat1})
for $X^{tx}$, we obtain: $$v_{i,n}(t,x)\leq v_{i,n+1}(t,x)\leq C_K (1+|x|^{\gamma}),$$ for some constants $C_K$ and $\gamma$ independent of $n$.
In order to complete the proof it is enough now to set $v_i(t,x):=\lim_{n\rightarrow \infty}v_{i,n}(t,x),\ \ (t,x)\in [0,T]\times \mathbb{R}^{k}$ since $Y^{i,n,tx}\nearrow Y^{i,tx}$ as $n \rightarrow \infty$. $\qed$
\medskip

We are now going to focus on the continuity of the functions $v_{1},...,v_{m}$.

\begin{thm}\label{continuity}
Assume that (H2), (H3) and (H4) are fulfilled. Then the functions $(v_{1},...,v_{m})(t,x):[0,T]\times \mathbb{R}^{k} \rightarrow \mathbb{R}$
 are continuous and solutions in viscosity sense of the system of variational inequalities with inter-connected obstacles (\ref{pdi}).
\end{thm}
$Proof$: Let us show that for every $i\in \mathcal{I}$, $v_{i}$ is continuous. First we have that
$v_{i}(t,x)=Y^{i,tx}_{t}$, then:
\begin{equation}\label{yitx-estimate}
\left|Y^{i,tx}_{t}-Y^{i,t^{'}x^{'}}_{t^{'}}\right|\leq
\left|Y^{i,tx}_{t}-Y^{i,n,tx}_{t}\right|+\left|Y^{i,n,tx}_{t}-Y^{i,n,t^{'}x^{'}}_{t^{'}}\right|+
\left|Y^{i,n,t^{'}x^{'}}_{t^{'}}-Y^{i,t^{'}x^{'}}_{t^{'}}\right|.
\end{equation}
Moreover,
\begin{eqnarray}\label{yintx-estimate}
\left|Y^{i,n,tx}_{t}-Y^{i,n,t^{'}x^{'}}_{t^{'}}\right| & \leq &\left|Y^{i,n,tx}_{t}-Y^{i,n,t^{'}x^{'}}_{t}\right|+\left|Y^{i,n,t^{'}x^{'}}_{t}-Y^{i,n,t^{'}x^{'}}_{t^{'}}\right| \nonumber \\
&& \nonumber \\
&\leq & \sup\limits_{0 \leq s \leq
T}\left|Y^{i,n,tx}_{s}-Y^{i,n,t^{'}x^{'}}_{s}\right| +
\left|Y^{i,n,t^{'}x^{'}}_{t}-Y^{i,n,t^{'}x^{'}}_{t^{'}}\right|.
\end{eqnarray}

Next we we will show the $L^{2}$--continuity of the value functions $(t,x)\rightarrow Y^{i,n,tx}$.\\
Recall that
\begin{equation*}
Y^{i,n,tx}_{s}=\limfunc{ess\
sup}_{u\in\mathcal{D}_{s}^{i,n}}E\left[\int_{s}^{T}\psi_{u_{r}}(r,X_{r}^{tx})dr-\sum_{j=1}^{n}g_{u_{\tau_{j-1}u_{\tau_{j}}}}
(\tau_{j},X_{\tau_{j}}^{tx})\mathbf{1}_{[\tau_{j}<T]}|\mathcal{F}_{s}\right],
\end{equation*}
where $\mathcal{D}_{s}^{i,n}=\{ u\in \mathcal{D} \ such \ that \
u_{0}=i, \tau_{1}\geq s, \ and \ \tau_{n+1}=T \}$. Therefore
\begin{eqnarray*}
\left|Y^{i,n,tx}_{s}-Y^{i,n,t^{'}x^{'}}_{s}\right|&=&\left|\limfunc{ess\sup}_{u\in\mathcal{D}_{s}^{i,n}}E\left[\int_{s}^{T}\psi_{u_{r}}(r,X_{r}^{tx})dr-\sum_{j=1}^{n}
g_{u_{\tau_{j-1}u_{\tau_{j}}}}(\tau_{j},X_{\tau_{j}}^{tx})\mathbf{1}_{[\tau_{j}<T]}|\mathcal{F}_{s}\right]\right.
\\ && \\
&&\left.-\limfunc{ess\
sup}_{u\in\mathcal{D}_{s}^{i,n}}E\left[\int_{s}^{T}\psi_{u_{r}}(r,X_{r}^{t^{'}x^{'}})dr-\sum_{j=1}^{n}g_{u_{\tau_{j-1}u_{\tau_{j}}}}
(\tau_{j},X_{\tau_{j}}^{t^{'}x^{'}})\mathbf{1}_{[\tau_{j}<T]}|\mathcal{F}_{s}\right]\right| \\ && \\
&\leq& \limfunc{ess\ sup}_{u\in\mathcal{D}_{s}^{i,n}}E\bigg[\int_{s}^{T}\left|\psi_{u_{r}}(r,X_{r}^{tx})-\psi_{u_{r}}(r,X_{r}^{t^{'}x^{'}})\right|dr \\ && \\
&&+\sum_{j=1}^{n}\left|\left(g_{u_{\tau_{j-1}u_{\tau_{j}}}}(\tau_{j},X_{\tau_{j}}^{tx})-g_{u_{\tau_{j-1}u_{\tau_{j}}}}(\tau_{j},X_{\tau_{j}}^{t^{'}x^{'}})\right)
\mathbf{1}_{[\tau_{j}<T]}\right||\mathcal{F}_{s}\bigg] \\ && \\
&\leq& E\left[\int_{0}^{T}\limfunc{max}_{i\in \mathcal{I}}\left|\psi_{i}(r,X_{r}^{tx})-\psi_{i}(r,X_{r}^{t^{'}x^{'}})\right|dr\right. \\ && \\
&&\left.+n\limfunc{max}_{l,j \in \mathcal{I};\ j\neq
l}\left\{\limfunc{sup}_{0\leq r\leq
T}\left|g_{lj}(r,X_{r}^{tx})-g_{lj}(r,X_{r}^{t^{'}x^{'}})\right|\right\}|\mathcal{F}_{s}\right].
\end{eqnarray*}

Now using Doob's Maximal Inequality and taking expectation, there
exists a constant $C\geq 0$ such that:
\begin{eqnarray}\label{esup}
&&E\left[\sup\limits_{0 \leq s \leq
T}\left|Y^{i,n,tx}_{s}-Y^{i,n,t^{'}x^{'}}_{s}\right|^{2}\right]\\
&\leq
& C
E\left[\int_{0}^{T}\limfunc{max}_{i\in\mathcal{I}}\left|\psi_{i}(r,X_{r}^{tx})-\psi_{i}(r,X_{r}^{t^{'}x^{'}})\right|^2
dr+n\limfunc{max}_{l,j \in \mathcal{I};\ j\neq
l}\left\{\limfunc{sup}_{0\leq r\leq
T}\left|g_{lj}(r,X_{r}^{tx})-g_{lj}(r,X_{r}^{t^{'}x^{'}})\right|\right\}
^2 \right]. \nonumber
\end{eqnarray}

In the right-hand side of (\ref{esup}) the first term converges to
$0$ as $(t',x')\rightarrow (t,x)$. Indeed, for any $\rho > 0$ it
holds true that:
\begin{eqnarray}
&&E\left[\int_{0}^{T}\limfunc{max}\limits_{i\in\mathcal{I}}\left|\psi_{i}(r,X_{r}^{tx})-\psi_{i}(r,X_{r}^{t^{'}x^{'}})\right|^{2}dr
\right]\nonumber\\
&\leq&
E\left[\int_{0}^{T}\limfunc{max}\limits_{i\in\mathcal{I}}\left|\psi_{i}(r,X_{r}^{tx})-\psi_{i}(r,X_{r}^{t^{'}x^{'}})\right|^{2}
\mathbf{1}_{[\mid X_{r}^{tx}\mid+\mid X_{r}^{t^{'}x^{'}}\mid \leq
\rho]}dr \right]\nonumber\\
&&+E\left[\int_{0}^{T}\limfunc{max}\limits_{i\in\mathcal{I}}\left|\psi_{i}(r,X_{r}^{tx})-\psi_{i}(r,X_{r}^{t^{'}x^{'}})\right|^{2}
\mathbf{1}_{[\mid X_{r}^{tx}\mid+\mid X_{r}^{t^{'}x^{'}}\mid >
\rho]}dr \right].\nonumber
\end{eqnarray} By the Lebesgue Dominated Convergence
Theorem, the continuity of $\psi_i$ and estimates (\ref{estimat2}),
the first term of the right-hand side of this inequality converges
to $0$ as $(t',x')\rightarrow (t,x)$.\\ The second term satisfies:
\begin{eqnarray}
&&E\left[\int_{0}^{T}\limfunc{max}\limits_{i\in\mathcal{I}}\left|\psi_{i}(r,X_{r}^{tx})-\psi_{i}(r,X_{r}^{t^{'}x^{'}})\right|^{2}
\mathbf{1}_{[\mid X_{r}^{tx}\mid+\mid X_{r}^{t^{'}x^{'}}\mid >
\rho]}dr \right]\nonumber
\\
&\leq& \left
\{E\left[\int_{0}^{T}\limfunc{max}\limits_{i\in\mathcal{I}}\left|\psi_{i}(r,X_{r}^{tx})-\psi_{i}(r,X_{r}^{t^{'}x^{'}})\right|^{4}dr
\right] \right\}^\frac{1}{2} \left \{E\left[\int_{0}^{T}
\mathbf{1}_{[\mid X_{r}^{tx}\mid+\mid X_{r}^{t^{'}x^{'}}\mid
> \rho]}dr \right]\right\}^\frac{1}{2}\nonumber
\\
 &\leq& \left
\{E\left[\int_{0}^{T}\limfunc{max}\limits_{i\in\mathcal{I}}\left|\psi_{i}(r,X_{r}^{tx})-\psi_{i}(r,X_{r}^{t^{'}x^{'}})\right|^{4}dr
\right] \right\}^\frac{1}{2} \left \{\rho^{-1}E\left[\int_{0}^{T}
(\mid X_{r}^{tx}\mid+\mid X_{r}^{t^{'}x^{'}}\mid )dr
\right]\right\}^\frac{1}{2}.\nonumber
\end{eqnarray}
 Using estimates
(\ref{estimat1}) and the polynomial growth of $\psi_i$, then we
have:
$$\begin{array}{ll}E\left[\int_{0}^{T}\limfunc{max}\limits_{i\in\mathcal{I}}\left|\psi_{i}(r,X_{r}^{tx})-\psi_{i}(r,X_{r}^{t^{'}x^{'}})\right|^{2}
\mathbf{1}_{[\mid X_{r}^{tx}\mid+\mid X_{r}^{t^{'}x^{'}}\mid >
\rho]}dr \right]\leq C(1+\mid x \mid^{\gamma}+\mid x'
\mid^{\gamma})\rho^{-\frac{1}{2}},
\end{array}$$
where $C$ and $\gamma$ are real constants which are bound to the
polynomial growth of $\psi_i$ and estimate (\ref{estimat1}).\\
As $\rho$ is arbitrary then making $\rho\rightarrow +\infty $ to
obtain that:

$$\lim \limits_{(t,x)\rightarrow(t',x')}E\left[\int_{0}^{T}\limfunc{max}\limits_{i\in\mathcal{I}}\left|\psi_{i}(r,X_{r}^{tx})-\psi_{i}(r,X_{r}^{t^{'}x^{'}})\right|^{2}
\mathbf{1}_{[\mid X_{r}^{tx}\mid+\mid X_{r}^{t^{'}x^{'}}\mid >
\rho]}dr \right] \rightarrow 0.$$ Thus the claim is proved.\\
In the same way we have:
$$\lim \limits_{(t,x)\rightarrow(t',x')} E\left[\limfunc{max}_{l,j \in \mathcal{I};\ j\neq
l}\left\{\limfunc{sup}_{0\leq r\leq
T}\left|g_{lj}(r,X_{r}^{tx})-g_{lj}(r,X_{r}^{t^{'}x^{'}})\right|\right\}
^2 \right]\rightarrow 0.$$ Then the right hand side of (\ref{esup})
converges to 0 as $(t^{'},x^{'})\rightarrow (t,x)$. Thus we obtain
that:
\begin{equation}\label{yitx-estimate2}
E\left[\sup\limits_{0 \leq s \leq
T}\left|Y^{i,n,tx}_{s}-Y^{i,n,t^{'}x^{'}}_{s}\right|^{2}\right]\rightarrow
0 \ as \ (t^{'},x^{'})\rightarrow (t,x).
\end{equation}
Then the function $(s,t,x)\rightarrow Y^{i,n,tx}_{s}$ is continuous
from $[0,T]^{2}\times \mathbb{R}^{k}$ into $L^{2}$, which is the
desired result.\\ Next, recall that from (\ref{yitx-estimate}) and
(\ref{yintx-estimate}) we have:
\begin{equation}\label{convergence}
\left|Y^{i,tx}_{t}-Y^{i,t^{'}x^{'}}_{t^{'}}\right|\leq
\left|Y^{i,tx}_{t}-Y^{i,n,tx}_{t}\right|+\sup\limits_{0 \leq s \leq
T}\left|Y^{i,n,tx}_{s}-Y^{i,n,t^{'}x^{'}}_{s}\right| +
\left|Y^{i,n,t^{'}x^{'}}_{t}-Y^{i,n,t^{'}x^{'}}_{t^{'}}\right|+\left|Y^{i,n,t^{'}x^{'}}_{t^{'}}-Y^{i,t^{'}x^{'}}_{t^{'}}\right|.
\end{equation}
We put $ n \rightarrow +\infty $ and using the fact that
$Y^{i,tx}_{t}$ is deterministic, $Y^{i,n,tx}$ converges to
$Y^{i,tx}$ locally uniformly and the continuity of $Y^{i,n,tx}_{t}$
in $t$, together with (\ref{yitx-estimate2}) we get that the right
hand side terms of (\ref{convergence}) converge to 0 as
$(t^{'},x^{'})\rightarrow (t,x)$. Therefore
$Y^{i,t^{'}x^{'}}_{t^{'}} \rightarrow Y^{i,tx}_{t}$ as
$(t^{'},x^{'})\rightarrow (t,x)$. Next since by Proposition
\ref{viscosity}, we have $\forall i\in \mathcal{I} \quad
Y^{i,tx}_{t}=v_i(t,x)$, then the deterministic functions
$v_{1},\ldots,v_m$ are continuous in $(t,x)$, moreover they are of
polynomial growth. Then taking into account Theorem \ref{yitx}
implies that $(v_{1},\ldots,v_m)$ is a viscosity solution for the
system of variational inequalities with inter-connected obstacles
(\ref{pdi}). The proof of Theorem \ref{continuity} is now complete.\
\qed

\section{Uniqueness of the solution of the system of Variational Inequalities}
In this section we are going to show the uniqueness of the viscosity solution of the system (\ref{pdi}). We first need the following lemma.
\begin{lem}\label{lemma-supersol}
Let $(v_{i})_{i=1,\ldots,m}$ be a supersolution of the system (\ref{pdi}), then for any $\gamma\geq 0$ there exists $\alpha>0$ such that for any $\lambda \geq \alpha$ and $\theta >0$, the m-uplet $(v_{i}(t,x)+\theta e^{-\lambda t}\mid x\mid^{2\gamma+2})_{i=1,\ldots,m}$ is a supersolution for (\ref{pdi}).
\end{lem}
$Proof. $ We assume w.l.o.g. that the functions
$(v_{i}(t,x))_{i=1,\ldots,m}$ are lsc. Let $ i\in \mathcal{I}$ be
fixed and let $\varphi \in \mathcal{C}^{1,2}$ be such that the
function $\varphi-(v_{i}+\theta e^{-\lambda t}\mid
x\mid^{2\gamma+2})$ has a local maximum in (t,x) which is equal to
0. Since $(v_{i}(t,x))_{i=1,\ldots,m}$ is a supersolution for
(\ref{pdi}), then we have: $\forall i\in \mathcal{I}$,
\begin{eqnarray}
&&\min\bigg\{ v_{i}(t,x)-\max\limits_{j\in \mathcal{I}^{-i}}\left\{-g_{ij}(t,x)+v_{j}(t,x)\right\},\nonumber\\&&\hspace{10mm}-\partial_{t}\Big(\varphi(t,x)-\theta e^{-\lambda t}\mid x\mid^{2\gamma+2}\Big)-\frac{1}{2}Tr\Big[\sigma.\sigma^{*}(t,x)D^{2}_{x}\Big(\varphi(t,x)-\theta
e^{-\lambda t}\mid x\mid^{2\gamma+2}\Big)\Big]\nonumber\\&&\hspace{10mm}-D_{x}\Big(\varphi(t,x)-\theta
e^{-\lambda t}\mid x\mid^{2\gamma+2}\Big).b(t,x)-\phi_{i}(t,x)\bigg\}
\geq 0.\nonumber
\end{eqnarray}
Hence
\begin{eqnarray}\label{supersol1}
&&(v_{i}(t,x)+\theta e^{-\lambda t}\mid
x\mid^{2\gamma+2})-\max\limits_{j\in
\mathcal{I}^{-i}}\Big(-g_{ij}(t,x)+(v_{j}(t,x)+\theta e^{-\lambda t}\mid
x\mid^{2\gamma+2})\Big) \nonumber %\\%&& \nonumber
\\&=&v_{i}(t,x)-\max\limits_{j\in
\mathcal{I}^{-i}}(-g_{ij}(t,x)+v_{j}(t,x)) \geq 0.
\end{eqnarray}
On the other hand:
\begin{equation*}
\begin{array}{l}
\left.-\partial_{t}\Big(\varphi(t,x)-\theta e^{-\lambda t}\mid x\mid^{2\gamma+2}\Big)-\frac{1}{2}Tr\Big[\sigma.\sigma^{*}(t,x)D^{2}_{x}\Big(\varphi(t,x)-\theta e^{-\lambda t}\mid x\mid^{2\gamma+2}\Big)\Big] \right.\\
\left.-D_{x}\Big(\varphi(t,x)-\theta e^{-\lambda t}\mid
x\mid^{2\gamma+2}\Big).b(t,x)-\phi_{i}(t,x)\right. \geq 0.
\end{array}
\end{equation*}
Thus
\begin{eqnarray}\label{supersol}
&&-\partial_{t}\varphi(t,x)-\frac{1}{2}Tr\Big[\sigma.\sigma^{*}(t,x)D^{2}_{x}\varphi(t,x)\Big]-D_{x}(\varphi(t,x)).b(t,x)-\phi_{i}(t,x) \nonumber\\
&\geq& \theta \lambda e^{-\lambda t}\mid
x\mid^{2\gamma+2}-\frac{1}{2}\theta e^{-\lambda t}
Tr\Big[\sigma.\sigma^{*}(t,x)D^{2}_{x}\mid x\mid^{2\gamma+2}\Big]-\theta
e^{-\lambda t}D_{x}(\mid x\mid^{2\gamma+2}).b(t,x).
\end{eqnarray}
Therefore taking into account the growth conditions on $b$ and $\sigma$ and setting $\theta>0$, there exists two positive constants $C_{1}$ and $C_{2}$ such that:\\
$\frac{1}{2}\theta e^{-\lambda t} Tr[\sigma.\sigma^{*}(t,x)D^{2}_{x}\mid x\mid^{2\gamma+2}]\leq C_{1} \mid x\mid^{2\gamma+2}$ and $D_{x}(\mid x\mid^{2\gamma+2}).b(t,x)\leq C_{2}\mid x\mid^{2\gamma+2}$. \\
Then by (\ref{supersol}) we get
\begin{equation}\label{supersol2}
\left.-\partial_{t}\varphi(t,x)-\frac{1}{2}Tr\Big[\sigma.\sigma^{*}(t,x)D^{2}_{x}\varphi(t,x)\Big]-D_{x}(\varphi(t,x)).b(t,x)-\phi_{i}(t,x) \right.\\
\left.\geq \theta (\lambda-\alpha)  e^{-\lambda t}\mid
x\mid^{2\gamma+2},\right.
\end{equation}
where $\alpha=C_{1}+C_{2}$.\\
Now since $\theta>0$, then we conclude that for $\lambda \geq \alpha$ the right hand side of (\ref{supersol2}) is non-negative.\\
Finally, noting that $i$ is arbitrary in $\mathcal{I}$ together with
(\ref{supersol1}), we obtain that \\$(v_{i}(t,x)+\theta e^{-\lambda
t}\mid x\mid^{2\gamma+2})_{i=1,\ldots,m}$ is a viscosity
supersolution for (\ref{pdi}).\qed \medskip

Now we give an equivalent form of the quasi-variational inequality
(\ref{pdi}). In this section, we consider the new function
$\Gamma_i$ given by the classical change of variable $\Gamma_i(t,x)
= \exp(t)v_i(t, x)$, for any $t\in[0,T ]$ and $x\in \mathbb{R}^{k}$.
Of course, the function $\Gamma_i$ is continuous and of polynomial
growth with respect to its second argument.\\ A second property is given
by the following proposition:
\begin{prop}
$v_i$ is a viscosity solution of (\ref{pdi}) if and only if
$\Gamma_i$ is a viscosity solution to the following
quasi-variational inequality in $[0,T [\times \mathbb{R}^{k}$,
\begin{equation} \label{sysvi10} \left\{
\begin{array}{ll}
\min \{\Gamma_i(t,x)- \max\limits_{j\in{\cal
I}^{-i}}(-e^{t}g_{ij}(t,x)+\Gamma_j(t,x)),\\
\qquad\qquad\qquad\qquad\qquad\qquad\Gamma_i(t,x)-\partial_t\Gamma_i(t,x)-
{\cal A}\Gamma_i(t,x)-e^{t}\psi_i(t,x)\}=0,\\
\Gamma_i(T,x)=e^{T}v_i(T,x)=0.\qed
\end{array}\right.
\end{equation}
\end{prop}
\medskip
\quad We are going now to address the question of uniqueness of the
viscosity solution of the system (\ref{pdi}). We notice that assumption (H4) will not be used in the proof of the uniqueness in the theorem below, and there are no restrictions on the negative switching costs.\\
We now give the main result of this section:
\begin{thm}\label{uniqueness}
Assume that (H2) and (H3)(i)--(ii) hold. Then the viscosity solution of the system of variational inequalities
with interconnected obstacles (\ref{pdi}) is unique in the space of
continuous functions on $[0,T]\times \mathbb{R}^{k}$ which satisfy a
polynomial growth condition, i.e., in the space
$$\begin{array}{l}{\cal C}:=\bigg\{\varphi: [0,T]\times \mathbb{R}^{k}\rightarrow
\mathbb{R}, \mbox{ continuous and for any }\\\qquad
\qquad\qquad(t,x), \, |\varphi(t,x)|\leq C(1+|x|^\gamma) \mbox{ for
some constants } C \mbox{ and }\gamma\bigg\}.\end{array}$$
\end{thm}
$Proof. $ We will show by contradiction that if $u_1,...,u_m$ and
$v_1,...,v_m$ are a subsolution and a supersolution respectively for
(\ref{sysvi10}) then for any $i=1,...,m$, $u_i\leq v_i$. Therefore
if we have two solutions of (\ref{sysvi10}) then they are obviously
equal. Next according to Lemma \ref{lemma-supersol}, it is enough to
show that for any $i \in \mathcal{I}$, we have:
$$\forall (t,x) \in [0,T]\times \mathbb{R}^{k}, \ \ u_{i}(t,x)\leq v_{i}(t,x)+\theta e^{-\lambda t}\mid x\mid^{2\gamma+2} \ ,$$
since in taking the limit as $\theta \rightarrow 0$ we get the
desired result. \\ So let us set  $w_{i}(t,x)=v_{i}(t,x)+\theta
 e^{-\lambda t}\mid x\mid^{2\gamma+2}$, $(t,x) \in [0,T]\times \mathbb{R}^{k}$.
  Next assume there exists a point $(\bar{t},\bar{x}) \in [0,T]\times \mathbb{R}^{k}$ such that for $i \in \mathcal{I}$: $\limfunc{max}_{i\in\mathcal{I}}(u_{i}(\bar{t},\bar{x})-w_{i}(\bar{t},\bar{x}))>0$. Then using the growth condition there exists $R>0$ such that:
$$\forall (t,x) \in [0,T]\times \mathbb{R}^{k} s.t. \mid x\mid \geq R, \ \ u_{i}(t,x)-w_{i}(t,x)<0 .$$
Since $u_{i}(T,x)=v_{i}(T,x)=0$, \ it implies that
\begin{eqnarray}\label{comp_uni}
0< \Max_{(t,x) \in [0,T]\times
\mathbb{R}^{k}}\limfunc{max}_{i\in\mathcal{I}}(u_{i}(t,x)-w_{i}(t,x))&=&
\Max_{(t,x) \in [0,T[\times
 B_R}\Max_{i\in\mathcal{I}}(u_{i}(t,x)-w_{i}(t,x))\nonumber\\
&=&\limfunc{max}_{i\in\mathcal{I}}(u_{i}(\widehat{t},\widehat{x})-w_{i}(\widehat{t},\widehat{x})),
\end{eqnarray}
where $B_R := \{x\in \mathbb{R}^{k}; |x|<R\}$ and $(\widehat{t},\widehat{x}) \in [0,T[\times B_R$.\\
Now let us define $\tilde{\mathcal{I}}$ as:
$$\tilde{\mathcal{I}}:=\bigg\{j \in \mathcal{I}, \ u_{j}(\widehat{t},\widehat{x})-w_{j}(\widehat{t},\widehat{x})=
\limfunc{max}_{i\in\mathcal{I}}(u_{i}(\widehat{t},\widehat{x})-w_{i}(\widehat{t},\widehat{x}))\bigg\}.$$
First note that $\tilde{\mathcal{I}}$ is not empty. Here $\gamma$ is
the growth exponent of the functions which w.l.o.g we assume integer
and $\geq 2$. Next for a small $\epsilon>0$ and $j \in
\tilde{\mathcal{I}}$, let us set for $(t,x,y) \in [0,T]\times B_R
\times B_R$,
$$\phi_{\epsilon}^{j}(t,x,y)=u_{j}(t,x)-w_{j}(t,y)-\varphi_{\epsilon}(t,x,y),$$
where,\ $$\varphi_{\epsilon}(t,x,y)=\frac{1}{2\epsilon}\mid
x-y\mid^{2\gamma}+\eta(\mid x-\widehat{x}\mid^{2\gamma+2}+ \mid
y-\widehat{x}\mid^{2\gamma+2}) + \beta(t-\widehat{t})^2$$ and
$\beta, \eta
>0$. Now let $(t_{\epsilon},x_{\epsilon},y_{\epsilon}) \in [0,T]\times \overline{B}_R \times \overline{B}_R$ be such that
:$$\phi^{j}_{\epsilon}(t_{\epsilon},x_{\epsilon},y_{\epsilon})=
\limfunc{max}_{(t,x,y) \in [0,T]\times \overline{B}_R\times
\overline{B}_R}\phi^{j}_{\epsilon}(t,x,y)$$ which exists since
$\phi^{j}_{\epsilon}$ is continuous. On the other hand, from
$2\phi^{j}_{\epsilon}(t_{\epsilon},x_{\epsilon},y_{\epsilon})\geq
\phi^{j}_{\epsilon}(t_{\epsilon},x_{\epsilon},x_{\epsilon})+\phi^{j}_{\epsilon}(t_{\epsilon},y_{\epsilon},y_{\epsilon})$,
we have
\begin{equation}
\frac{1}{2\epsilon}|x_{\epsilon} -y_{\epsilon}|^{2\gamma} \leq
(u_{j}(t_{\epsilon},x_{\epsilon})-u_{j}(t_{\epsilon},y_{\epsilon}))+(w_{j}(t_{\epsilon},x_{\epsilon})-w_{j}(t_{\epsilon},y_{\epsilon})),
\end{equation}
and consequently $\frac{1}{2\epsilon}|x_{\epsilon}
-y_{\epsilon}|^{2\gamma}$ is bounded, and as $\epsilon\rightarrow
0$, $|x_{\epsilon} -y_{\epsilon}|\rightarrow 0$. Since $u_{j}$
and $w_{j}$ are uniformly continuous on $[0,T]\times \overline{B}_R$, then $\frac{1}{2\epsilon}|x_{\epsilon} -y_{\epsilon}|^{2\gamma}\rightarrow 0$ as $\epsilon\rightarrow 0.$\\
Since
 \begin{equation}\label{uj-wj}
u_{j}(\widehat{t},\widehat{x})-w_{j}(\widehat{t},\widehat{x})=\phi^{j}_{\epsilon}(\widehat{t},\widehat{x},\widehat{x})\leq
\phi^{j}_{\epsilon}(t_{\epsilon},x_{\epsilon},y_{\epsilon}) \leq
u_{j}(t_{\epsilon},x_{\epsilon})-w_{j}(t_{\epsilon},y_{\epsilon}),
\end{equation}
it follows as $\epsilon\rightarrow 0$ and the continuity of $u$ and
$w$ that, up to a subsequence,
 \begin{equation}\label{subsequence}
 (t_\epsilon,x_\epsilon,y_\epsilon)\rightarrow (\widehat{t},\widehat{x},\widehat{x}).
 \end{equation}
We now claim that for some  $k \in \tilde{\mathcal{I}}$ we have:
$$u_{k}(\widehat{t},\widehat{x})>\limfunc{max}_{j \in \mathcal{I}^{-k}}(u_{j}(\widehat{t},\widehat{x})-e^{\widehat{t}}g_{kj}(\widehat{t},\widehat{x})).$$
Indeed if for any $k \in \tilde{\mathcal{I}}$ we have:
$$u_{k}(\widehat{t},\widehat{x})\leq \limfunc{max}_{j \in \mathcal{I}^{-k}}(u_{j}(\widehat{t},\widehat{x})-e^{\widehat{t}}g_{kj}(\widehat{t},\widehat{x})),$$
then there exists $j \in \mathcal{I}^{-k}$ such that:
$$u_{k}(\widehat{t},\widehat{x})- u_{j}(\widehat{t},\widehat{x})\leq -e^{\widehat{t}}g_{kj}(\widehat{t},\widehat{x}).$$
From the supersolution property of $w_{j}$, we have
$$w_{k}(\widehat{t},\widehat{x})\geq \limfunc{max}_{j \in \mathcal{I}^{-k}}(w_{j}(\widehat{t},\widehat{x})-e^{\widehat{t}}g_{kj}(\widehat{t},\widehat{x})).$$
Then
$$w_{k}(\widehat{t},\widehat{x})-w_{j}(\widehat{t},\widehat{x})\geq -e^{\widehat{t}}g_{kj}(\widehat{t},\widehat{x}).$$
It follows that
$$u_{k}(\widehat{t},\widehat{x})- u_{j}(\widehat{t},\widehat{x})\leq -e^{\widehat{t}}g_{kj}(\widehat{t},\widehat{x})\leq w_{k}(\widehat{t},\widehat{x})-w_{j}
(\widehat{t},\widehat{x}).$$ Hence
$$u_{k}(\widehat{t},\widehat{x})-w_{k}(\widehat{t},\widehat{x})\leq -e^{\widehat{t}}g_{kj}(\widehat{t},\widehat{x})+u_{j}(\widehat{t},\widehat{x})-w_{k}
(\widehat{t},\widehat{x})\leq
u_{j}(\widehat{t},\widehat{x})-w_{j}(\widehat{t},\widehat{x}).$$ But
since $k\in \tilde{\mathcal{I}}$ then
$$u_{k}(\widehat{t},\widehat{x})-w_{k}(\widehat{t},\widehat{x})= u_{j}(\widehat{t},\widehat{x})-w_{j}(\widehat{t},\widehat{x})=
 -e^{\widehat{t}}g_{kj}(\widehat{t},\widehat{x})+u_{j}(\widehat{t},\widehat{x})-w_{k}(\widehat{t},\widehat{x}),$$
which implies that j belongs also to $\tilde{\mathcal{I}}$ and
$$u_{k}(\widehat{t},\widehat{x})-u_{j}(\widehat{t},\widehat{x})=-e^{\widehat{t}}g_{kj}(\widehat{t},\widehat{x}).$$
Repeating this procedure as many times as necessary and since
$\tilde{\mathcal{I}}$ is finite we get the existence of a loop of
indices $i_{1},...,i_{p},i_{p+1}$ of $\tilde{\mathcal{I}}$ such that
$i_{p+1}=i_{1}$ and
$$g_{i_{1},i_{2}}(\widehat{t},\widehat{x})+...+g_{i_{p},i_{p+1}}(\widehat{t},\widehat{x})=0.$$
But this contradicts the assumption (\ref{nofreeloop}), whence the claim holds. \\
To proceed let us consider $k\in \tilde{\mathcal{I}}$ such that:
\begin{equation}
u_{k}(\widehat{t},\widehat{x})>\limfunc{max}_{j \in
\mathcal{I}^{-k}}(u_{j}(\widehat{t},\widehat{x})-g_{kj}(\widehat{t},\widehat{x})).
\end{equation}
By the continuity of $u_{j}$ and $g_{ij}$ and since
$(t_{\epsilon},x_{\epsilon},u_{k}(t_{\epsilon},x_{\epsilon}))\rightarrow_{\epsilon}
(\widehat{t},\widehat{x},u_{k}(\widehat{t},\widehat{x}))$ then for
$\epsilon$ small enough we have:
\begin{equation}\label{uktnxn}
u_{k}(t_{\epsilon},x_{\epsilon})>\limfunc{max}_{j \in
\mathcal{I}^{-k}}(u_{j}(t_{\epsilon},x_{\epsilon})-g_{kj}(t_{\epsilon},x_{\epsilon})).
\end{equation}
Next we have : \begin{equation}\left\{
\begin{array}{lllll}\label{derive}
D_{t}\varphi_{\epsilon}(t,x,y)=2\beta(t-\widehat{t}),\\
D_{x}\varphi_{\epsilon}(t,x,y)=
\frac{\gamma}{\epsilon}(x-y)|x-y|^{2\gamma-2} +\eta(2\gamma + 2)
(x-\widehat{x})|x-\widehat{x}|^{2\gamma}, \\
D_{y}\varphi_{\epsilon}(t,x,y)=
-\frac{\gamma}{\epsilon}(x-y)|x-y|^{2\gamma-2} +
\eta(2\gamma + 2)(y-\widehat{x})|y-\widehat{x}|^{2\gamma},\\\\
B(t,x,y)=D_{x,y}^{2}\varphi_{\epsilon}(t,x,y)=\frac{1}{\epsilon}
\begin{pmatrix}
a_1(x,y)&-a_1(x,y) \\
-a_1(x,y)&a_1(x,y)
\end{pmatrix}+ \begin{pmatrix}
a_2(x)&0 \\
0&a_2(y)
\end{pmatrix}, \\\\
\mbox{ with } a_1(x,y)=\gamma|x-y|^{2\gamma-2}I+\gamma(2\gamma -2)(x-y)(x-y)^* |x-y|^{2\gamma-4} \mbox{ and }\\
a_2(x)=\eta(2\gamma +
2)|x-\widehat{x}|^{2\gamma}I+2\eta\gamma(2\gamma +
2)(x-\widehat{x})(x-\widehat{x})^* |x-\widehat{x}|^{2\gamma-2 }.
\end{array}
\right. \end{equation} Taking into account (\ref{uktnxn}) then
applying the result by Crandall et al. (Theorem 8.3 in  \cite{[CIL]})
to the function $$
u_{k}(t,x)-w_{k}(t,y)-\varphi_{\epsilon}(t,x,y) $$ at the
point $(t_\epsilon,x_\epsilon,y_\epsilon)$, for any $\epsilon_1 >0$,
we can find $c,d\in \mathbb{R}$ and $X,Y \in S_k$, such that:

\begin{equation} \label{lemmeishii} \left\{
\begin{array}{lllll}
(c,\frac{\gamma}{\epsilon}(x_\epsilon-y_\epsilon)|x_\epsilon-y_\epsilon|^{2\gamma-2}
+\eta(2\gamma +
2)(x_\epsilon-\widehat{x})|x_\epsilon-\widehat{x}|^{2\gamma},X)
\in J^{2,+}(u_{k}(t_\epsilon,x_\epsilon)),\\
(-d,\frac{\gamma}{\epsilon}(x_\epsilon-y_\epsilon)|x_\epsilon-y_\epsilon|^{2\gamma-2}
-\eta(2\gamma +
2)(y_\epsilon-\widehat{x})|y_\epsilon-\widehat{x}|^{2\gamma },Y)\in
J^{2,-}w_{k}(t_\epsilon,y_\epsilon)),\\
c+d=D_{t}\varphi_{\epsilon}(t_\epsilon,x_\epsilon,y_\epsilon)=2\beta(t_\epsilon-\widehat{t}) \mbox{ and finally }\\
-(\frac{1}{\epsilon_1}+||B(t_\epsilon,x_\epsilon,y_\epsilon)||)I\leq
\begin{pmatrix}
X&0 \\
0&-Y
\end{pmatrix}\leq B(t_\epsilon,x_\epsilon,y_\epsilon)+\epsilon_1 B(t_\epsilon,x_\epsilon,y_\epsilon)^2.
\end{array}
\right. \end{equation} Taking now into account (\ref{uktnxn}), and
the definition of viscosity solution, we get:
$$\begin{array}{l}-c+u_{k}(t_\epsilon,x_\epsilon)-\frac{1}{2}Tr[\sigma^{*}(t_\epsilon,x_\epsilon)X\sigma(t_\epsilon,x_\epsilon)]
-\langle\frac{\gamma}{\epsilon}(x_\epsilon-y_\epsilon)|x_\epsilon-y_\epsilon|^{2\gamma-2}
+\\\qquad\qquad\qquad\qquad\qquad\eta(2\gamma +
2)(x_\epsilon-\widehat{x})|x_\epsilon-\widehat{x}|^{2\gamma
},b(t_\epsilon,x_\epsilon)\rangle-e^{t_\epsilon}\psi_{k}(t_\epsilon,x_\epsilon)\leq
0 \mbox{ and
}\\d+w_{k}(t_\epsilon,y_\epsilon)-\frac{1}{2}Tr[\sigma^{*}(t_\epsilon,y_\epsilon)Y\sigma(t_\epsilon,y_\epsilon)]
-\langle\frac{\gamma}{\epsilon}(x_\epsilon-y_\epsilon)|x_\epsilon-y_\epsilon|^{2\gamma-2}
-\\\qquad\qquad\qquad\qquad\qquad\eta(2\gamma +
2)(y_\epsilon-\widehat{x})|y_\epsilon-\widehat{x}|^{2\gamma},b(t_\epsilon,y_\epsilon)\rangle-e^{t_\epsilon}\psi_{k}(t_\epsilon,y_\epsilon)\geq
0\end{array}$$ which implies that:
\begin{equation}
\begin{array}{llll}
\label{viscder} -c-d
+u_{k}(t_\epsilon,x_\epsilon)-w_{k}(t_\epsilon,y_\epsilon)&\leq
\frac{1}{2}Tr[\sigma^{*}(t_\epsilon,x_\epsilon)X\sigma
(t_\epsilon,x_\epsilon)-\sigma^{*}(t_\epsilon,y_\epsilon)Y\sigma(t_\epsilon,y_\epsilon)]\\
&\qquad +
\langle\frac{\gamma}{\epsilon}(x_\epsilon-y_\epsilon)|x_\epsilon-y_\epsilon|^{2\gamma-2},b(t_\epsilon,x_\epsilon)
-b(t_\epsilon,y_\epsilon)\rangle\\&\qquad+\langle \eta(2\gamma +
2)(x_\epsilon-\widehat{x})|x_\epsilon-\widehat{x}|^{2\gamma
},b(t_\epsilon,x_\epsilon)\rangle\\&\qquad +\langle \eta(2\gamma +
2)(y_\epsilon-\widehat{x})|y_\epsilon-\widehat{x}|^{2\gamma
},b(t_\epsilon,y_\epsilon)\rangle
\\&\qquad+e^{t_\epsilon}\psi_{k}(t_\epsilon,x_\epsilon)-e^{t_\epsilon}\psi_{k}(t_\epsilon,y_\epsilon).
\end{array}
\end{equation}
But from (\ref{derive}) there exist two constants $C$ and $C_1$ such
that:
$$||a_1(x_\epsilon,y_\epsilon)||\leq C|x_\epsilon - y_\epsilon|^{2\gamma -2} \mbox{ and }(||a_2(x_\epsilon)||\vee ||a_2(y_\epsilon)||)\leq C_1\eta .$$
As
$$B= B(t_\epsilon,x_\epsilon,y_\epsilon)= \frac{1}{\epsilon}
\begin{pmatrix}
a_1(x_{\epsilon},y_{\epsilon})&-a_1(x_{\epsilon},y_{\epsilon}) \\
-a_1(x_{\epsilon},y_{\epsilon})&a_1(x_{\epsilon},y_{\epsilon})
\end{pmatrix}+ \begin{pmatrix}
a_2(x_\epsilon)&0 \\
0&a_2(y_\epsilon)
\end{pmatrix},$$
then
$$B\leq \frac{C}{\epsilon}|x_\epsilon - y_\epsilon|^{2\gamma -2}
\begin{pmatrix}
I&-I \\
-I&I
\end{pmatrix}+ C_1 \eta I.$$
It follows that:
\begin{equation}
B+\epsilon_1 B^2 \leq C(\frac{1}{\epsilon}|x_\epsilon -
y_\epsilon|^{2\gamma -2}+ \frac{\epsilon_1}{\epsilon^2}|x_\epsilon -
y_\epsilon|^{4\gamma -4})\begin{pmatrix}
I&-I \\
-I&I
\end{pmatrix}+ C_1 \eta I,
\end{equation}
where $C$ and $C_1$ which hereafter may change from line to line.
Choosing now $\epsilon_1=\epsilon$, yields the relation
\begin{equation}
\label{ineg_matreciel} B+\epsilon_1 B^2 \leq
\frac{C}{\epsilon}(|x_\epsilon - y_\epsilon|^{2\gamma
-2}+|x_\epsilon - y_\epsilon|^{4\gamma -4})\begin{pmatrix}
I&-I \\
-I&I
\end{pmatrix}+ C_1 \eta I.
\end{equation}
Now, from (\ref{b}), (\ref{lemmeishii}) and (\ref{ineg_matreciel}) we get:
$$\frac{1}{2}Tr[\sigma^{*}(t_\epsilon,x_\epsilon)X\sigma(t_\epsilon,x_\epsilon)-\sigma^{*}(t_\epsilon,y_\epsilon)
Y\sigma(t_\epsilon,y_\epsilon)]\leq \frac{C}{\epsilon}(|x_\epsilon -
y_\epsilon|^{2\gamma}+|x_\epsilon - y_\epsilon|^{4\gamma -2}) +C_1
\eta (1+|x_\epsilon|^2+|y_\epsilon|^2).$$ Next $$
\langle\frac{\gamma}{\epsilon}(x_\epsilon-y_\epsilon)|x_\epsilon-y_\epsilon|^{2\gamma-2},b(t_\epsilon,x_\epsilon)-b(t_\epsilon,y_\epsilon)\rangle
\leq \frac{C^2}{\epsilon}|x_\epsilon - y_\epsilon|^{2\gamma},$$ and
finally,
\begin{eqnarray}
&&\langle \eta(2\gamma + 2)(x_\epsilon-\widehat{x})|x_\epsilon-\widehat{x}|^{2\gamma},b(t_\epsilon,x_\epsilon)\rangle +
\langle \eta(2\gamma +
2)(y_\epsilon-\widehat{x})|y_\epsilon-\widehat{x}|^{2\gamma
},b(t_\epsilon,y_\epsilon)\rangle\nonumber
\\
&\leq& C\eta(1+|x_\epsilon||x_\epsilon-\widehat{x}|^{2\gamma +
1}+|y_\epsilon||y_\epsilon-\widehat{x}|^{2\gamma + 1}).\nonumber
\end{eqnarray}
So that by plugging into (\ref{viscder}) we obtain:
\begin{eqnarray*}
&&-2\beta(t_\epsilon-\widehat{t})+ u_{k}(t_\epsilon,x_\epsilon)-w_{k}(t_\epsilon,y_\epsilon)\\
&\leq&
\frac{C}{\epsilon}(|x_\epsilon - y_\epsilon|^{2\gamma}+|x_\epsilon -
y_\epsilon|^{4\gamma -2}) +C_1 \eta
(1+|x_\epsilon|^2+|y_\epsilon|^2) +\frac{C^2}{\epsilon}|x_\epsilon -
y_\epsilon|^{2\gamma}\\&&+C\eta
(1+|x_\epsilon||x_\epsilon-\widehat{x}|^{2\gamma +
1}+|y_\epsilon||y_\epsilon-\widehat{x}|^{2\gamma +
1})+
e^{t_\epsilon}\psi_{k}(t_\epsilon,x_\epsilon)-e^{t_\epsilon}\psi_{k}(t_\epsilon,y_\epsilon).
\end{eqnarray*}
By sending $\epsilon\rightarrow0$, $\eta \rightarrow0$ and taking
into account the continuity of $\psi_{k}$ and $\gamma \geq 2$, we
obtain:
$$u_k(\widehat{t},\widehat{x})-w_k(\widehat{t},\widehat{x})\leq 0,$$
 which contradicts (\ref{comp_uni}). The proof of Theorem
\ref{uniqueness} is now complete. \qed


\begin{thebibliography}{99}
\small
\renewcommand{\baselinestretch}{0.3}
\bibitem{[BE]} E. Bayraktar and M. Egami, On the One-Dimensional Optimal Switching Problem, {Mathematics of Operations Research,} 2010, 35 (1), 140-159.

\bibitem {[BOU]} {\small B. Bouchard, A stochastic target formulation for optimal
switching problems in finite horizon, {Stochastics,} 81 (2009), pp.
171--197. }

\bibitem{[BO1]} {\small K. A. Brekke and B. {\O }ksendal, Optimal switching in
an economic activity under uncertainty. {\ SIAM J. Control Optim,}
(32) (1994), pp. 1021-1036. }


\bibitem{[BS]}  {\small M. J. Brennan and E. S. Schwartz, Evaluating natural
resource investments, {J. Business} 58 (1985), pp. 135--137. }


\bibitem{[CK]} {\small J. Cvitanic and I. Karatzas, Backward SDEs with
reflection and Dynkin games, {\ Annals of Probability} 24 (4)
(1996), pp. 2024--2056. }


\bibitem{[CL]} {\small R. Carmona and M. Ludkovski, Pricing asset scheduling
flexibility using optimal switching, {Appl. Math. Finance}, 15
(2008), pp. 405--447. }

\bibitem{[CIL]} {\small M. Crandall, H. Ishii, H and P.L. Lions, User's guide
to viscosity solutions of second order partial differential
equations, Bull. Amer. Math. Soc, 27 (1992), 1--67. }

\bibitem{[DM]} {\small C. Dellacherie and P. A. Meyer, Probabilit\'{e}s et Potentiel, V-VIII, Hermann, Paris, 1980.}

\bibitem{[DX]} {\small S. J. Deng and Z. Xia, Pricing and Hedging Electric
Supply Contracts: A Case with Tolling Agreements, preprint, {Georgia
Institute of Technology}, Atlanta, 2005. }

\bibitem{[D]} {\small A. Dixit, Entry and exit decisions under uncertainty, {
J. Political Economy}, 97 (1989), pp. 620--638. }

\bibitem{[DP]} {\small A. Dixit and R. S. Pindyck, Investment Under
Uncertainty, Princeton University Press, Princeton, NJ, 1994. }

\bibitem{[DH]} {\small B. Djehiche and S. Hamad\`{e}ne, On a finite horizon
starting and stopping problem with risk of abandonment, {Int. J.
Theor. Appl. Finance,} 12 (2009), pp. 523--543.}

\bibitem{[DHP]} {\small B. Djehiche, S. Hamad\`{e}ne, A. Popier, A finite
horizon optimal multiple switching problem, {SIAM J. Control Optim}.
48 (4) (2009) 2751--2770. }


\bibitem{[DZ]} {\small K. Duckworth and M. Zervos, A problem of stochastic
impulse control with discretionary stopping, in Proceedings of the
39th IEEE Conference on Decision and Control, {IEEE Control Systems
Society}, Piscataway, NJ, 2000, pp. 222--227. }

\bibitem{[DZ2]} {\small K. Duckworth and M. Zervos, A model for investment
decisions with switching costs, {Ann. Appl. Probab}, 11 (2001), pp.
239--260. }

\bibitem{[E1]} {\small B. El Asri, Stochastic Optimal Multi-Modes Switching with a
Viscosity Solution Approach, {Stochastic Processes
and their Applications}, 123 (2013), pp. 579--602.}

\bibitem{[E2]} {\small B. El Asri, Deterministic minimax impulse control in finite horizon: the viscosity solution approach,
 {ESAIM: Control, Optimisation and Calculus of Variations}, 19 (2013), pp. 63--77.}

\bibitem{[EH]} {\small B. El Asri and S. Hamad\`{e}ne, The finite horizon
optimal multi-modes switching problem: The viscosity solution
approach,{\ Appl. Math. Optim}, 60 (2009), pp. 213--235. }

\bibitem{[Elka]} {\small N. El Karoui, Les aspects probabilistes du contr\^ole
stochastique, in Ecole d'\'et\'e de Probabilit\'es de Saint-Flour,
Lecture Notes in Math. 876, Springer-Verlag, New York, 1980. }

\bibitem{[EKal]} {\small N. El Karoui, C. Kapoudjian, E. Pardoux, S. Peng, and
M. C. Quenez, Reflected solutions of backward SDEs and related
obstacle problems for PDEs, {Ann. Probab.,} 25 (1997), pp. 702--737.
}

\bibitem{[HJ]} {\small S. Hamad\`{e}ne and M. Jeanblanc, On the starting and
stopping problem: Application in reversible investments, Math. Oper.
Res., 32 (2007), pp. 182--192. }

\bibitem{[HM]} {\small S. Hamad\`ene and M.A. Morlais, Viscosity Solutions of Systems of PDEs with
Interconnected Obstacles and Multi--Modes Switching Problem, {Applied Mathematics and Optimization},
67 (2013), pp. 163--196.}

\bibitem{[HZ]} {\small S. Hamad\`{e}ne and J. Zhang, Switching problem and
related system of reflected backward SDEs, Stochastic Processes and
their Applications, 120 (2010) 403--426. }

\bibitem{[HT]} {\small Y. Hu and S. Tang, Multi-dimensional BSDE with oblique
reflection and optimal switching, Probab. Theory Related Fields
(2009) doi:10.1007/s00440-009-0202-1. }

\bibitem{[LP]} {\small V. Ly Vath and H. Pham, Explicit solution to an optimal
switching problem in the two-regime case, {SIAM J. Control Optim.,}
46 (2007), pp. 395--426. }

\bibitem{[RM]} {\small T. S. Knudsen, B. Meister, and M. Zervos, Valuation of
investments in real assets with implications for the stock prices,
SIAM J. Control Optim., 36 (1998), pp. 2082--2102. }

\bibitem{[PTW]} {\small A. Porchet, N. Touzi, and X. Warin, Valuation of
power plants by utility indifference and numerical computation,
Math. Methods Oper. Res., 70 (2009), pp. 47--75. }

\bibitem{[RY]} {\small D. Revuz and M. Yor, Continuous Martingales and
Brownian Motion, Springer-Verlag, Berlin, 1991. }

\bibitem{[shi]} {\small H. Shirakawa, Evaluation of investment opportunity
under entry and exit decisions. \textit{S${\bar u}$rikaisekikenky${\bar u}$%
sho K${\bar o}$ky${\bar u}$roku} (987) (1997), pp. 107--124. }

\bibitem{[TY]} {\small S. Tang and J. Yong, Finite horizon stochastic optimal
switching and impulse controls with a viscosity solution approach,
Stoch. Stoch. Rep., 45 (1993), pp. 145--176. }

\bibitem{tri1} {\small L. Trigeorgis, Real options and interactions with
financial flexibility, Financial Management, 22 (1993), pp.
202--224. }

\bibitem{[tri]} {\small L. Trigeorgis, Real Options: Managerial Flexibility
and Strategy in Resource Allocation, MIT Press, Cambridge, MA, 1996.
}

\bibitem{[dz]} {\small M. Zervos, A problem of sequential entry and exit
decisions combined with discretionary stopping, SIAM J. Control
Optim., 42 (2003), pp. 397--421. }


%\bibitem{ht} Hu, Y., Tang, S. (2007):
%Multi-dimensional BSDE with Oblique Reflection and Optimal
%Switching. {\it Preprint Universit\'e de Rennes 1, France}
%
%\bibitem Ishii, H., and Koike, K.:viscosity solutions of a system of Nonlinear second order PDE's arising in switching games, Functional , 34, 143-155, 1991.
%
%
%\bibitem{[LP]} Ly Vath, V. and Pham, H. (2007): Explicit solution to an optimal switching problem in the two-regime case.
%SIAM Journal on Control and Optimization, pp. 395-426.
%
%\bibitem{[RM]} Knudsen, T. S., Meister, B. and Zervos, M. (1998):
%Valuation of investments in real assets with implications for the stock prices. {\it SIAM J. Control Optim.}
%(36), pp. 2082-2102.
%
%\bibitem {[PH]} Pham, H.: On the smooth-fit property for one-dimensional optimal switching problem.
%
%\bibitem{porchet} Porchet, A., Touzi, N., Warin, X. (2006):
%Valuation of a power plant under production constraints. {\it
%Preprints of the 10th Annual Conference in Real Options, NYC, USA,
%June, pp. 14-17,
%http://wwww.realoptions.org/abstracts/abstracts06.html}
%
%\bibitem{porchet2} Porchet, A., Touzi, N., Warin, X. (2007): Valuation of a Power Plant Under Production
%Constraints and Market Incompleteness, {\it to appear in Management
%Science (2008)}
%
%\bibitem{[RY]} Revuz, D and Yor, M. (1991): Continuous Martingales and
%Brownian Motion. {\it Springer Verlag, Berlin.}
%
%\bibitem{shi} Shirakawa, H. (1997): Evaluation of investment
%opportunity under entry and exit decisions. {\it S${\bar u}$rikaisekikenky${\bar u}$sho K${\bar o}$ky${\bar
%u}$roku} (987), pp. 107-124.
%
%\bibitem{[TY]} Tang, S. and Yong, J. (1993): Finite horizon stochastic optimal switching and impulse
%controls with a viscosity solution approach.  Stoch. and Stoch. Reports, 45, 145-176.
%
%\bibitem{tri1} Trigeorgis, L. (1993): Real options and interactions
%with financial flexibility. {\it Financial Management} (22), pp. 202-224.
%
%\bibitem{tri} Trigeorgis, L. (1996): Real Options: Managerial
%Flexibility and Startegy in Resource Allocation. {\it MIT Press.}
%
%\bibitem{[dz]} Zervos, M. (2003): A Problem of Sequential Enty and Exit
%Decisions Combined with Discretionary Stopping. {\it SIAM J. Control Optim.} 42 (2), pp. 397-421.

\end{thebibliography}
\end{document}